\documentclass[12pt,oneside]{article}
\newcommand{\Path}{}
\newcommand{\figs}{}
\usepackage{ \Path Paper_style}
\usepackage{ \Path Keywords}
\usepackage{ \Path Environments}
\usepackage{tikz-cd}
\usepackage{amsmath}
\usepackage{tabularx} \newcolumntype{L}{>{\raggedright\arraybackslash}X}
\usepackage{lineno}

\DeclareMathOperator{\Feature}{Feature}
\title{Data-driven discovery of quasiperiodically driven dynamics}
\author{Suddhasattwa Das, Shakib Mustavee, Shaurya Agarwal }

\begin{document}
	\maketitle
	\begin{abstract} The analysis of a timeseries can provide many new perspectives if it is accompanied by the assumption that the timeseries is generated from an underlying dynamical system. For example, statistical properties of the data can be related to measure theoretic aspects of the dynamics, and one can try to recreate the dynamics itself. The underlying dynamics could represent a natural phenomenon or a physical system, where the timeseries represents a sequence of measurements. In this paper, we present a completely data-driven framework to identify and model quasiperiodically driven dynamical systems (Q.P.D.) from the timeseries it generates. Q.P.D. are a special class of systems that are driven by a periodic source with multiple base frequencies. Such systems abound in nature, e.g., astronomy and traffic flow. Our framework reconstructs the dynamics into two components - the driving quasiperiodic source with generating frequencies; and the driven nonlinear dynamics. We make a combined use of a kernel-based harmonic analysis, kernel-based interpolation technique, and Koopman operator theory. Our framework provides accurate reconstructions and frequency identification for three real-world case studies.
	\end{abstract}
	
	
	
	\section{Introduction} \label{sec:intro}
	
	In this article we present a broad framework for analyzing timeseries generated by dynamical systems, to reconstruct the underlying dynamics as well as extract salient features of the dynamics. 
	A dynamical system can be described most generally as a space $\Omega$, along with a map / transformation $F:\Omega \to \Omega$. Depending on whether $\Omega$ is identified as a topological space, probability space, or a manifold, the study of the dynamics is labeled as topological, ergodic, or differential dynamical systems theory. Any point $z_0\in \Omega$ can be interpreted as a state of the system, and $F(z_0)$ is to be interpreted as the immediate next state of the system. Repeated applications of the dynamics $F$ leads to an orbit 
	\[ z_n := F^n(z_0) , \quad  n=1,2,3,\ \ldots , \]
	which is a sequence of points in $\Omega$. Such a general setup can be made into a hypothesis for timeseries analysis. We state this formally : 
	
	\begin{Assumption} \label{A:data}
		There is an unknown continuous function/observation $Y: \Omega \to \real^{\alpha}$, which is possibly a low-dimensional / partial observation of the dynamics. The data available for processing is the sequence of $\alpha$-dimensional data point $\SetDef{ y_n := Y(z_n) }{ n=0,1,2,\ldots }$, where $(z_n)$ is a trajectory of the dynamics under $F$.
	\end{Assumption}
	
	Thus timeseries analysis can be interpreted as an indirect study of the dynamics, with the only information being available being the timeseries $\braces{y_n}_{n=1,2,\ldots}$. This is the motivation for the field of \emph{data-driven} discovery of dynamical systems \cite[e.g.]{garcia2023physics, zhang2023application, DGJ_compactV_2018}. Due to the immense variety of dynamical behavior, such techniques are usually limited in scope. Effective techniques are designed towards more specific goals, such as control of parameterized systems \cite[e.g.]{Rahmani2023fractional, wu2023driver, DSSY2017_QR}, or to extract salient features such as patterns and correlations \citep[e.g.][]{afzali2023resonances, MustaveeEtAl_covid_2021, GiannakisDas_tracers_2019}. Our goal is to provide a general and completely non-parametric method, which only based on some assumptions on the dynamics itself. These assumptions are utilized to sharpen the data-driven algorithms.
	
	\paragraph{Skew-products} Our focus is on a class of dynamical systems which we call \emph{quasiperiodically driven dynamics}, stated formally as the map
	\begin{equation} \begin{split} \label{eqn:def:quasi_driven}
			F : \mathbb{T}^d \times \mathcal{M} \to \mathbb{T}^d \times \mathcal{M}, \quad 
			\paran{ \begin{array}{c} \theta_{n+1} \\ x_{n+1} \end{array} } = 
			F \paran{ \begin{array}{c} \theta_{n} \\ x_{n} \end{array} } =
			\paran{ \begin{array}{c} \theta_n + \rho \bmod \mathbb{T}^d \\ g( \theta_n, x_n ) \end{array} } .
		\end{split},\end{equation}
	Here the variable $\theta$ is an angular coordinate on a $d$-dimensional torus $\mathbb{T}^d$, $x$ is a point in some abstract or unknown manifold $\mathcal{M}$, and $g : \mathbb{T}^d \times \mathcal{M} \to \mathcal{M}$ is some nonlinear function. The coordinate labelled $\theta$ represents the \emph{phase} of a driving quasiperiodic rotation \cite{DasJim2017_SuperC, Dioph_Herman_1979}, and the vector $\rho$ is called the \emph{rotation vector} \citep[e.g.][]{Herman1, Arnold1965}. The coordinate $\rho$ represents the angular increments at each step for each of the coordinates of $\theta$. Thus \eqref{eqn:def:quasi_driven} is a \emph{one-way coupled} or \emph{skew-product} dynamical system on the space $\Omega := \mathbb{T}^d \times \mathcal{M}$. This model of a dynamical system captures the essence of many physical systems which are influenced or driven by an external, quasiperiodic phenomena. Some examples are mechanical systems \cite[e.g.]{kundu2023nonlinear}, astronomical data \cite{DSSY2017_QR},  
	climate data \cite{VautardGhil89, SlawinskaGiannakis16}, and physical flows on periodic domains \cite{GiannakisDas_tracers_2019, froyland2010coherent, FroylandEtAl14}. This autonomous periodic phenomenon could be seasonal, weekly or diurnal cycles, or geographic location. If the underlying system arises from a continuous time system by taking samples at intervals $\Delta t$, then $\rho = \Delta t \omega$ for some angular frequency vector $\omega$. This model can be summarized as the formal assumption :
	
	\begin{Assumption} \label{A:basic}
		There is a dynamical system of the form \eqref{eqn:def:quasi_driven} for some $\rho\in \mathbb{T}^d$, $m$-dimensional $C^2$ manifold $\mathcal{M}$, and a $C^2$ function $g : \mathbb{T}^d \times \mathcal{M} \to \mathcal{M}$. This dynamics has an invariant Borel probability measure $\mu$ with compact support $X\subseteq \mathbb{T}^d \times \mathcal{M}$.
	\end{Assumption}
	
	The structural assumption of \eqref{eqn:def:quasi_driven} has been shown to hold in great generality for dynamical systems with mixed spectrum \citep[see][Corr 7]{Das2023Koop_susp}. All dynamical systems which can be studied through experiments or measurements must have bounded trajectories. And any dynamical system with bounded trajectories must have at least one invariant probability measure \citep[e.g.][]{BrinStuck, KatokHassel1997}. We next utilize the components of Assumption \ref{A:basic} for a closer look at \eqref{eqn:def:quasi_driven}.
	
	\paragraph{Driven dynamics} Let $\nu$ be the push forward of $\mu$ onto $\mathcal{M}$, defined for every open subset $A$ of $\mathcal{M}$ as
	\[ \nu(A) := \mu \SetDef{ (\theta, x) }{ \theta\in \mathbb{T}^d, \, x\in A } = \mu \paran{ \mathbb{T}^d \times A } . \]
	%
	Now consider the following averages
	\begin{equation} \label{eqn:def:gper}
		g_{per}( \theta ) := \int g \paran{ \theta, x} d\nu(x), \quad \forall \theta\in \mathbb{T}^d.
	\end{equation}
	As a result, one can rewrite $g$ as
	\begin{equation} \label{eqn:odf3}
		g(\theta, x) = g_{per} (\theta) + g_{chaos} \paran{ \theta,x }, \quad \forall \theta\in \mathbb{T}^d, \, \forall x\in \mathcal{M}.
	\end{equation}
	The coordinate $\theta$ represents a periodically changing \emph{phase} of an autonomous dynamical system $\theta_{n+1} = \theta_n + \rho \bmod{ 2 \pi}$. The coordinate $x$ is to be interpreted as a set of variables, driven by $\theta$ and also simultaneously by its own current value. By virtue of \eqref{eqn:def:gper}, the component $g_{chaos}$ has zero-mean with respect to $\nu$, for every $\theta$. The system will be said to have \emph{constant sensitivity} if the following holds :
	
	\begin{Assumption}[Constant sensitivity] \label{A:cnstnt_snstv}
		The partial derivative of the function $g$ (from \eqref{eqn:def:quasi_driven} ) with respect to $x$ is independent of $\theta$. In other words
		\[ \frac{\partial}{ \partial \theta } \frac{\partial}{ \partial x } g \equiv 0 . \]
	\end{Assumption}

	Assumption~\ref{A:cnstnt_snstv} implies that $ \frac{\partial}{\partial x} g_{chaos}$ is a function of $x$ alone, so $g_{chaos}$ must take the form 
	\[g_{chaos}(\theta,x) = a(\theta) + \tilde{g}(x) .\]
	However \eqref{eqn:def:gper} and \eqref{eqn:odf3} imply that $\int g_{chaos}(\theta,x) d\nu(x)$ is zero for every $\theta$, which implies that $a(\theta) \equiv 0$. Thus $g_{chaos}$ may be interpreted as a function of $x$ alone, and \eqref{eqn:def:quasi_driven} simplifies into the following form :
	\begin{equation}\begin{split} \label{eqn:def:quasi_driven_II}
			\theta_{n+1} &= \theta_n + \rho \bmod 2\pi \\
			x_{n+1} &= g_{per}( \theta_n) + g_{chaos}(x_n)
	\end{split}\end{equation}
	Equation~\eqref{eqn:def:quasi_driven_II} interprets the action of $g$ on $x$ as a zero-mean function $g_{chaos}$ of $x$, along with a mean value dependent on the phase $\theta$ via the function $g_{per}$. The dynamics in the variable $x$ takes the format of additively forced dynamical systems, studied in various contexts \citep[e.g.][]{foster2006optimal, zhu2015tipping, DasJim17_chaos}.
	
	\paragraph{Goal} Our aim is to create a data-driven method that formulates a model in the format of \eqref{eqn:def:quasi_driven_II}. The iterations of this model should produce the same time series $\braces{y_n}_{n=1,2,\ldots}$ used in its construction. The task is made harder by the fact that none of the state variables $x,\theta$ are necessarily being observed. The function $Y$ $\braces{y_n}_{n=1,2,\ldots}$ is also unknown. To summarize, all the functions, spaces, and dimension $d$ described in Assumption~\ref{A:basic} and \ref{A:cnstnt_snstv} are unknown, and the only information available is the timeseries described in Assumption~\ref{A:data}. The specific objectives become : (i) finding the \emph{quasiperiodicity dimension} $d$ and the rotation vector $\omega$; and the functions (ii) $g_{per}$ and (iii) $g_{chaos}$.
	
\begin{figure}[!htbp]\center
		\includegraphics[width=.95\linewidth, height = 0.7 \textheight, keepaspectratio]{\figs 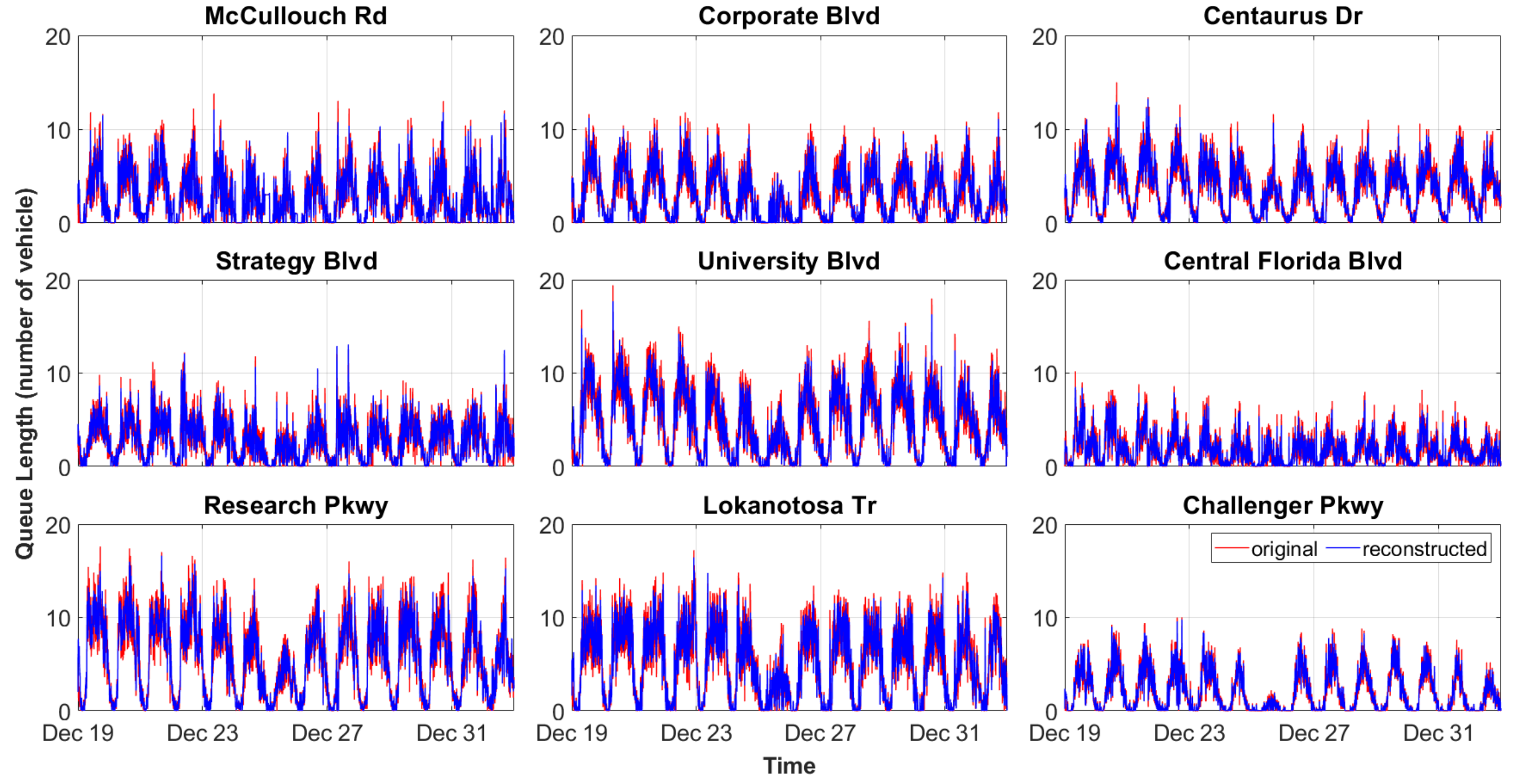}
		\caption{ Traffic intersection dynamics as a quasiperiodically driven system.  The red curves in each of the nine panels are traffic queue-length data collected at nine different intersections on the Alafaya corridor in Florida, USA. These curves carry the effects of stochastic, chaotic as well as periodic components of the traffic system. The blue curve is the output obtained from simulating a numerically constructed dynamical system described in \eqref{eqn:reconstruct}.  The close match between the two curves, in terms of their periodicities, fluctuations, and absolute error, highlights the theoretical and numerical contributions of our work. We begin with an axiomatical definition of a class of dynamics called \emph{quasiperiodically driven} dynamics, via Equations \eqref{eqn:def:quasi_driven}, \eqref{eqn:def:quasi_driven_II}, and formal Assumptions \ref{A:basic} and \ref{A:cnstnt_snstv}. We develop kernel-based learning techniques in Section~\ref{sec:data}, and show that in a purely data-driven setting described in Assumption \ref{A:data}, the results of our techniques converge to the true dynamical system. Our techniques rely heavily on ergodic theory and RKHS theory, described in Sections \ref{sec:eigen} and \ref{sec:kernel}, respectively. See Section \ref{sec:examples} for a full description of the traffic-intersection example, as well as two other physical systems. Also see Figures \ref{fig:408_weekly} and \ref{fig:atrial} for a similar analysis of other real-world dynamical systems as quasiperiodically driven dynamics. }
		\label{fig:Alafaya_reconstruct}
	\end{figure}

\begin{figure}[!htbp]\center
		\includegraphics[width=.8\linewidth]{\figs 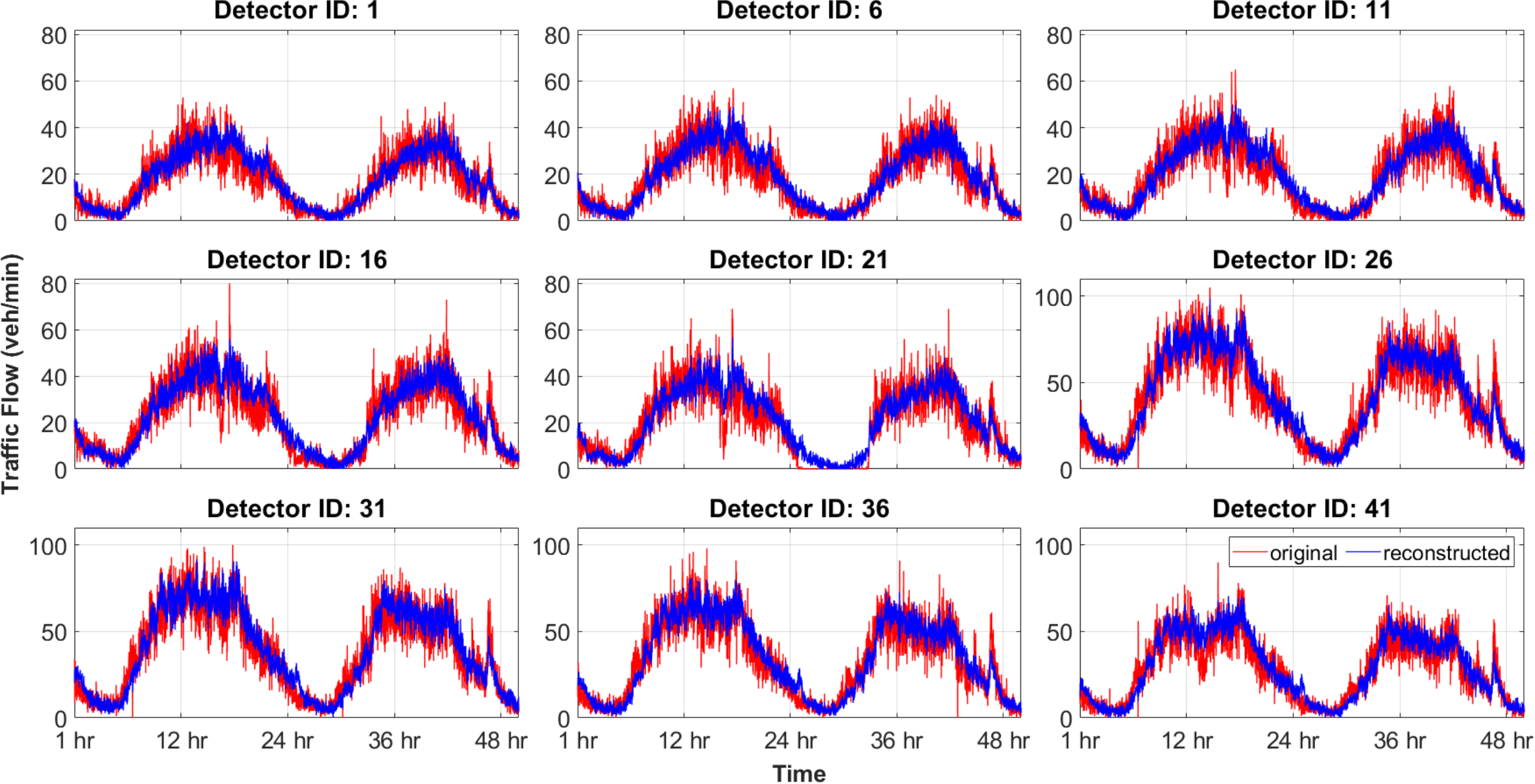}
		\caption{ Freeway traffic flow as a quasiperiodically driven system. The analysis here is similar to that in Figure \ref{fig:Alafaya_reconstruct}, but for traffic flow data collected by 47 detectors along the SR 408 freeway in Florida, USA. }
		\label{fig:408_weekly}
	\end{figure}

\begin{figure}[!htbp]\center
		\includegraphics[width=.95\linewidth]{\figs 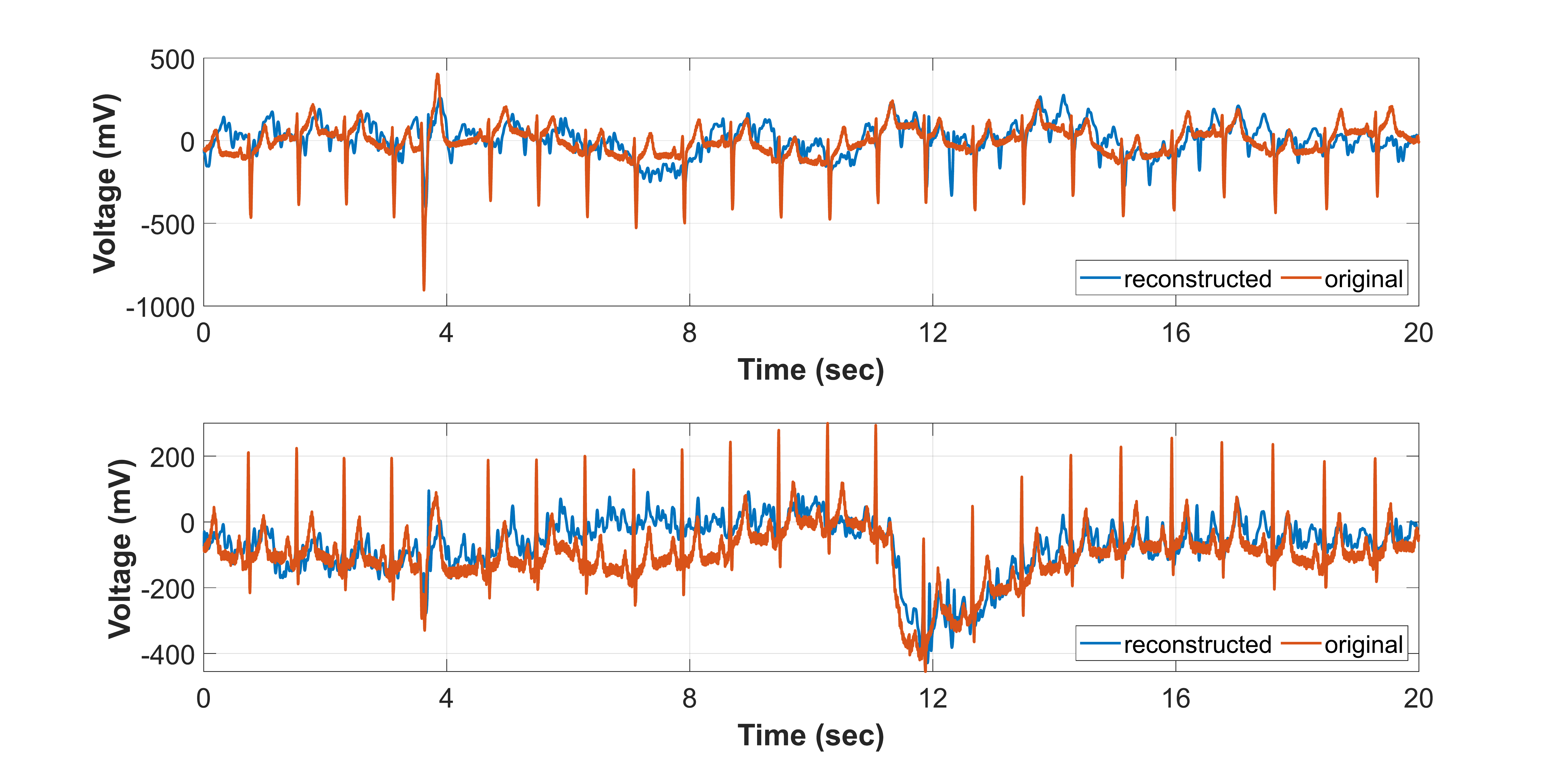}
		\caption{ Human cardiac system as a quasiperiodically driven system. The analysis here is similar to that in Figure \ref{fig:Alafaya_reconstruct}, but for atrial fibrillation (top) and atrial flutter (bottom) data obtained from patients. Unlike the traffic systems in the previous two examples, the heart has non-smooth dynamics, driven by periodic but sudden firing of neurons. We selected this example as a test for the limits of our learning algorithm. The reconstruction error is higher but correctly captures the periodicities and spiking events.
		}
		\label{fig:atrial}
	\end{figure}

	\paragraph{Approach} Our techniques rely on the Koopman operator theoretic formulation of the dynamics and are implemented using kernel-based techniques for learning and frequency analysis \cite{DasGiannakis_RKHS_2018}. See Figures~\ref{fig:Alafaya_reconstruct}, \ref{fig:408_weekly}, and \ref{fig:atrial} for illustrations of the results of applying our methods to various real-world systems. The numerical techniques we use are based on principles from ergodic theory \cite{DasGiannakis_RKHS_2018} and kernel-based learning theory \citep[e.g.]{PaulsenRaghupathi2016, Paulsen2016}. As a result, our reconstruction provably converges (see Theorem \ref{thm:2}) to the true dynamics, in a general situation precisely defined by a set of formal Assumptions. 
	

	\paragraph{Challenges and contributions} We have axiomatically defined a class of dynamics called \emph{quasiperiodically driven} dynamics, via Equations \eqref{eqn:def:quasi_driven}, \eqref{eqn:def:quasi_driven_II}, and formal Assumptions \ref{A:basic} and \ref{A:cnstnt_snstv}.  There are many challenges to accomplishing the goal of identifying the true frequency and reconstructing the original dynamics. Firstly, we show in Section~\ref{sec:eigen} that a dynamical system either has only one eigenfrequency (=0) or infinitely many. Moreover, if $d>1$ in \eqref{eqn:def:quasi_driven}, then the eigenfrequencies are dense on the real line. This makes numerically identifying and separating these eigenfrequencies challenging. Secondly, the presence of a mixing / chaotic component makes the traditional frequency analysis techniques such as dynamic mode decomposition \citep[e.g.][]{WilliamsEtAl15, KordaEtAl2018} or Fourier analysis \citep[e.g.][]{Wiener_Wintner_1941, Katznelson2004} unreliable. Mixing is the property of decay of correlations, see \cite{Nadkarni} or \citep[][Sec 3]{DasGiannakis_delay_2019} for more details. 
	To overcome these challenges, we obtain guarantees of convergence to the true discrete spectrum of the system using the proposed RKHS-based filtering technique (Algorithm~\ref{algo:RKHS_kernel}) described in Section~\ref{sec:data}. We further demonstrate through three extensive case studies (in Section \ref{sec:examples}) that in a purely data-driven setting described in Assumption \ref{A:data}, the results of our technique converge to the true dynamical system.
	Thirdly, we divide the task of learning / reconstructing the function $g$ from \eqref{eqn:def:quasi_driven} into its two components $g_{per}$ and $g_{chaos}$ from \eqref{eqn:def:quasi_driven_II}. A leaning approach oblivious to the internal structure described in \eqref{eqn:def:quasi_driven_II} could lead to a reconstructed dynamical system with widely divergent trajectories. 
	See Table~\ref{tab:compare} for a comparison of our methods with other techniques, particularly DMD-based techniques.
	
	
	\paragraph{Outline} A key consideration for us is \emph{quasiperiodicity}.  In Section~\ref{sec:eigen}, we discuss what it means and also interpret the significance of the $\theta$ coordinates. A key component of our method is the use of kernel integral operators and the theory of reproducing kernel Hilbert spaces. We discuss these concepts and the relevant techniques in Section~\ref{sec:kernel}. The actual data-driven implementation of the theory is described in Section~\ref{sec:data}, where we operate under Assumption \ref{A:data}. Finally, we use these techniques to analyze three real-world systems in Section~\ref{sec:examples}. Figures \ref{fig:Alafaya_reconstruct}, \ref{fig:408_weekly} and \ref{fig:atrial} reveal the results of these applications. The various symbols and notations we use are summarized in Table \ref{tab:param}.
	
	\begin{table}[!tbp]
		\caption{Mathematical Notations}
		\begin{tabular} { p{0.1\textwidth}   p{0.7\textwidth} }
			Notations & Description  \\ 
			\hline
			$F$ &  quasiperiodically driven dynamical system\\
			$d$ & dimension (degrees of freedom) of quasiperiodicity of $F$\\
			$\mathbb{T}^d$ & $d$-dimensional torus of driving quasiperiodic system\\
			$\mathcal{M}$ & manifold on which $x$ lies\\
			$\Omega$ & phase space $\mathbb{T}^d \times \mathcal{M}$\\
			$\theta$ & any point on ${T}^d$ representing the phase of the driving quasiperiodic system \\
			$\omega$ & angular frequency vector related to the rotation of the torus\\
			$\Delta t$ & sampling interval\\
			$\rho$ & rotation vector of the torus\\
			$x$ & driven states of the dynamical system \\
			$g$ & nonlinear function on $x$ and $\theta$ representing the driven dynamics \\
			$g_{per}$ & periodic component of $g$\\
			$g_{chaos}$ & chaotic component of $g$\\
			$U$ & the Koopman operator\\
			$\zeta$ & an eigenfunction of Koopman operator $U$\\
			$k$ & a kernel on space $M$\\
			$K$ & integral operator associated with kernel $k$\\
			$\alpha$ & dimension of the system measurements\\
			$y_n$ & sequence of $\alpha$-dimensional observation or data points\\ 
			$N$ & number of data points \\
			$(\theta_n, x_n)$ & trajectory of the dynamical system under $F$\\
			$\epsilon$ & Gaussian kernel bandwidth\\
			$ P_\epsilon$ & a compact, symmetric operator on $L^2(\mu)$ \\
			$Q$ & number of delay coordinates\\
			$L$ & number of eigenfunctions computed for the kernel integral matrix\\
			$L_0$ & number of eigenfunctions used to compute RKHS norm \\
			$\epsilon_1$ & threshold for the first filtration of candidate eigenfrequencies \\
			$\epsilon_2$ & threshold for the second filtration of candidate eigenfrequencies\\
			\hline
		\end{tabular}
		\label{tab:param1}
	\end{table}
	
	\section{Koopman operator and its spectrum} \label{sec:eigen}
	
	The Koopman operator $U$ converts the original nonlinear dynamics on a finite-dimensional phase-space into a linear dynamics on an infinite dimensional vector space. It is a time-shift operator, operating on functions instead of points on the phase space. Given any function $\phi:\Omega \to \real$, $U\phi$ is another function defined as
	\begin{equation} \label{eqn:def:Koop}
		(U\phi)(z) := \phi\left( F z\right) , \quad \forall z\in\Omega ,
	\end{equation}
	%
	%
	where $F$ is the underlying dynamical system \eqref{eqn:def:quasi_driven}. $\phi$ can be interpreted as a measurement or observation on the phase space $\Omega$, and $U\phi$ is the evolution / transformation of this measurement with the dynamics. Note that by virtue of \eqref{eqn:def:Koop} the correspondence $\phi \mapsto U\phi$ is linear. Thus the Koopman operator converts any nonlinear dynamical system into a linear map. This allows various tools from operator theory / functional analysis to be brought into the study of dynamics.  The properties of $U$ depend on the choice of vector/function space. Some common choices of function spaces are $C(\Omega)$ the space of conitnuous functions, or $C^r(\Omega)$, the space of $r$-times differentiable functions on $\Omega$. We shall use the Hilbert space $L^2(\mu)$, the space of square-integrable functions with respect to an invariant measure $\mu$ of the dynamics. In this space, the Koopman operator is a unitary operator \citep[see][Sec 2]{Nadkarni}, a property which makes its spectrum have some desirable properties and be numerically accessible \cite{DasGiannakis_delay_2019, DasGiannakis_RKHS_2018, DGJ_compactV_2018}.
	
	\paragraph{Koopman eigenfrequencies} Eigenvalue and eigenfunction pairs are one of the first attributes of an operator that are studied. For $U$, these carry a lot of significance. Since $U$ is unitary, its spectrum must lie on the unit circle of the complex plane. The eigenvalues of $U$ correspond to the point spectrum, and any eigenfunction $\zeta$ has a corresponding eigenvalue of the form $e^{\iota \omega}$ for some $\omega\in \real$. $\omega$ is called the Koopman eigenfrequency corresponding to $\zeta$. We thus have
	\begin{equation} \label{eqn:Koop_eig}
		(U^n \zeta)(z) \stackrel{\mbox{by def.}}{=} \zeta( F^n z) = e^{\iota \omega n} \zeta(z) , \quad \forall n\in\num .
	\end{equation}
	Equation~\ref{eqn:Koop_eig} reveals that the time-evolution of Koopman eigenfunctions is highly tractable, it is equivalent to multiplication by $e^{\iota \omega n}$ as a function of time $n$. Since $e^{\iota \omega}$ lies on the unit circle, the magnitude of $U^n \zeta$ does not change with $n$. As a result, the prediction formula for the evolution of an eigenfunction is not only simple (e.g. \eqref{eqn:Koop_eig}), it is also stable to initial approximation errors. $U$ always has the constant functions as eigenfunctions with eigenfrequency $0$. In general $U$ may or may not have other eigenfrequencies. For the special structure \eqref{eqn:def:quasi_driven} that we assume, the eigenfunctions of the driving system $\theta \mapsto \theta + \rho$ provide eigenfunctions for the dynamics under $F$. Suppose $\zeta$ is a Koopman eigenfunction for the driving system. Then we have :
	\[ \zeta(\theta_{n+k}) = e^{\iota\omega k} \zeta(\theta_n), \quad \forall k,n\in\num . \]
	Then define $\bar{\zeta}(x, \theta) := \zeta(\theta)$. Then note that
	\[ \bar{\zeta}( \theta_{n+k}, x_{n+k}) = \zeta(\theta_{n+k}) = e^{\iota\omega k} \zeta(\theta_n) = e^{\iota\omega k} \bar\zeta(\theta_n, x_n), \quad \forall k,n\in\num . \]
	Thus every eigenfunction for the driving system also leads to an eigenfunction for the entire system. 
	We can now give an alternate definition of the integer $d$ in \eqref{eqn:def:quasi_driven}, in terms of \emph{generating frequencies}. 
	
	\paragraph{Koopman eigenfunctions as phase} An important realization for us is that we may assume without loss of generality that $\abs{\zeta}\equiv 1$. This is because under the additional assumption of ergodicity \citep[see][]{DasGiannakis_delay_2019} we can assume that $\zeta$ is non-zero almost everywhere. Then if we set
	\[ \tilde{\zeta}(z) := \frac{\zeta(z)}{ \abs{\zeta(z)} } , \]
	it can be easily verified that $\tilde{\zeta}$ is also a Koopman eigenfunction with same frequency. Henceforth, we shall write a Koopman eigenfunction as a map %
	\[ \zeta : \Omega \to S^1 , \]
	where $S^1$ is the unit circle in the complex plane $\cmplx$. Thus the value a Koopman eigenfunction assigns to a point can be interpreted as the phase of the point. According to \eqref{eqn:Koop_eig}, this phase rotates uniformly with speed $\omega$. Koopman eigenfunctions thus reveal the rotational / (quasi)-periodic factors embedded in chaotic dynamics. This quasiperiodic component may not be evident from the state-space equations or from the measurement. However, the presence of such eigenfunctions strongly affects the outcome of various data-driven procedures  \citep[e.g.][]{DasGiannakis_delay_2019, DasJim2017_SuperC,  DSSY_Mes_QuasiP_2016}.
	
	This interpretation of eigenfunctions as phase helps reveal torus dynamics embedded within the system. Suppose we have a Koopman eigenfunction $\zeta:\Omega\to \cmplx$ with eigenfrequency $\omega$. Then by \eqref{eqn:Koop_eig}
	\[ \zeta( F^n x) = e^{\iota \omega n} \zeta(x) . \]
	We had discussed before that we can assume without loss of generality that $\abs{\zeta} \equiv 1$. This leads to the factored dynamics :
	\[\begin{tikzcd}
		\Omega \arrow{d}{\zeta} \arrow{r}{F} & \Omega  \arrow{d}{\zeta} \\
		S^1 \arrow{r}{R_\omega} & S^1
	\end{tikzcd}\]
	where $R_\omega$ is the uniform rotation by angular speed $\omega$ on the unit circle. Therefore if the eigenfrequency is nonzero, the values of $\zeta$ trace out the unit circle on the complex plane. The diagram above also indicates that the evolution of a Koopman eigenfunction is a dynamics of its own. For this reason, we shall use the terms ``Koopman mode", ``Koopman eigenmode" interchangeably with the term ``Koopman eigenfunction". Thus a Koopman mode reveals a circle rotation factored into the original dynamics. In fact, if we combine $m$ Koopman modes, we get a rotation on a $m$-dimensional torus :
	\begin{equation} \label{eqn:jne3}
		\begin{tikzcd}
			\Omega \arrow{d}[swap]{\zeta_1, \ldots, \zeta_m} \arrow{r}{F} & \Omega  \arrow{d}{\zeta_1, \ldots, \zeta_m} \\
			\mathbb{T}^m \arrow{r}{R_{\vec\omega}} & \mathbb{T}^m
		\end{tikzcd}
	\end{equation}
	If these eigenfrequencies are independent, then any orbit in $\Omega$ has a dense image in $\mathbb{T}^m$ under $\zeta_1, \ldots, \zeta_m$. Since the image of $\paran{\zeta_1, \ldots, \zeta_m}$ is closed, this makes them a surjective map. 
	Thus Koopman eigenfunctions reveal embedded toral dynamics of the same dimension as the quasiperiodicity dimension.
	
	\paragraph{Generating frequencies} The collection of eigenfunctions and (eigen)-frequencies have an algebraic structure to them. For any two  frequencies $\omega_1, \omega_2$, and integers $a,b$, $a\omega_1 + b\omega_2$ is also a frequency. This is because if $z_1, z_2$ are their corresponding eigenfunctions, then
	\[ z_1( F (\theta, x)) = e^{\iota \omega_1 } z_1(\theta, x), \quad z_2( F x) = e^{\iota \omega_2} z_2(\theta, x).  \]
	As a result
	\[ (z_1^a z_2^b) ( F (\theta, x)) =  z_1^a( F (\theta, x))  z_2^b( F (\theta, x)) = e^{\iota(a \omega_1 + b\omega_2)} (z_1^a z_2^b) (\theta, x) . \]
	Thus integer linear combinations of frequencies are again frequencies, and products of eigenfunctions are again eigenfrequencies. This makes the eigenfrequencies a module over the ring of integers. In particular, if the system has at least one nonzero frequency, then it has all harmonics of that frequency and thus infinitely many frequencies. A collection of eigenfrequencies is said to be independent if no integer linear combination of them is an integer.
	If the system has two independent frequencies, then all its frequencies are together dense on the real line. 
	
	A collection of frequencies will be called a \emph{basis} or \emph{generating} set of eigenfrequencies if they are independent and all frequencies of the system can be generated by taking integer linear combinations of frequencies from this set. There is no unique choice of a basis, but all bases will have the same dimension $d$, called the quasiperiodicity dimension $d$. In finite-dimensional manifolds such as $\Omega$, the number $d$ is usually observed to be finite \citep[e.g.][Sec 3]{DasGiannakis_delay_2019}, a fact that is supported by the fact that the factor map in \eqref{eqn:jne3} is surjective. Using results from Lie group theory \citep[][Thm 6]{Das2023Lie}, it can be shown that if a generating set of Koopman eigenfunctions are smooth, then the factor map in \eqref{eqn:jne3} becomes a submersion of manifolds. We show in Section~\ref{sec:data} how one can avoid the task of finding such a generating set, the quasiperiodicity dimension $d$, as well as rotation vector $\rho$, and still obtain a data-driven reconstruction.
	
	
	This completes our discussion on the quasiperiodic structure of the dynamics \eqref{eqn:def:quasi_driven}. We next discuss some techniques from Functional Analysis for reconstructing the quasiperiodic component and its complement.
	
	\section{Kernels and integral operators} \label{sec:kernel}
	
	A kernel is a function $ k : M\times M\to \real $ on some space $M$. The quantity $k(x,y)$ is a measure of similarity, closeness, or distance between two points $x,y\in M$. Kernel-based methods have been used very effectively to obtain geometric information of the underlying space $M$. This information has been used to study various related structures such as statistical manifolds  \cite{DasDimitEnik2020}, geometric information \cite{DimitrisBerry_SEC_2018, BerryHarlim2016, BerrySauer2017}, and dynamical information such as tracer flows \cite{GiannakisDas_tracers_2019}, Lyapunov functions \cite{GieslHafstein2015}, stable/unstable foliations \cite{BerryEtAl2013}, Koopman spectrum  \cite{Giannakis2015, DasGiannakis_delay_2019, DasGiannakis_RKHS_2018} and more generally the spectral measure \cite{DGJ_compactV_2018}. The techniques in this paper are based on \cite{DasGiannakis_RKHS_2018}. We shall use the \emph{Gaussian kernel}
	\[ \kappa_\epsilon(x,y) := \exp\left( -\frac{1}{\epsilon} d(x,y)^2 \right), \]
	where $\epsilon$ is called the \emph{bandwidth parameter}, and $d(\cdot, \cdot)$ is some notion of metric or distance on the space. Note that $k_\epsilon(x,y)=1$ iff $x=y$. It decays exponentially from $1$ as $y$ moves away from $x$. If $\epsilon$ is decreased, then the function decays more sharply. For our purposes, $M = \Omega$, the phase-space of the dynamics. However, since we are working under the data-driven assumption~\ref{A:data}, $\Omega$ will be assumed to be unknown and we need an indirect access to $\Omega$. This is done through an embedding described below.
	
	\paragraph{Delay-coordinates} The data sequence described in Assumption~\ref{A:data} is obtained through an observation $Y$. However, $Y$ may not faithfully replicate $\Omega$, i.e., $Y$ may not be a one-to-one map and its values may not correspond to unique states in $\Omega$. An easy solution to this problem is the method of delay coordinates \citep[e.g.][]{SauerEtAl1991, Sauer1992_delay, BerryDas_learning_2022}, in which the dynamics is embedded in higher dimensional space $\real^{\alpha(Q+1)}$, where $Q$ is called the number of delays. The delay coordinated version of the map $Y$ is the map
	\[ Y^{(Q)} : \Omega \to \real^{\alpha(Q+1)}, \quad Y^{(Q)}(\omega) := \left( Y(\omega), Y( F^1 \omega), \ldots, Y(F^Q\omega) \right) . \]
	Thus the delay coordinated version of each point $y_n$ is
	\[ y_n \leftrightarrow y_N^{(Q)} := \left( y_n, y_{n+1}, \ldots, y_{n+Q} \right) . \]
	The main point of using delay coordinates, as explained in \cite{SauerEtAl1991}, is that for a typical observation map $Y$, if $Q$ is large enough, then $Y^{(Q)} : \Omega \to \real^{\alpha(Q+1)}$ is an embedding / one-to-one map. There are several heuristic algorithms to determine a $Q$ which would be sufficiently large \citep[e.g.][]{BuzugPfister1992, BuzugPfister1992_pre, SauerYorke1993_embed, Aguirre1995_delays}. We proceed with the assumption that the $Q$ chosen is large enough. We can then use the Gaussian shape function to implicitly obtain a kernel $k_\epsilon : \Omega \times \Omega \to \real$ as follows 
	\begin{equation} \label{eqn:def:Gauss_ker} 
		k_\epsilon( z, z') := \kappa_\epsilon\left( Y^{(Q)} (z), Y^{(Q)} (z') \right) = \exp\left( -\frac{1}{\epsilon} \norm{ Y^{(Q)}(z) - Y^{(Q)} (z') }^2  \right). 
	\end{equation}
	
	Even if the two states $z, z'$ are unknown, the left-hand side in \eqref{eqn:def:Gauss_ker} can be computed since the right-hand side only uses the observation map $Y$. When using the Gaussian kernel directly with finite data, one can run into problems of undersampling or non-uniform density \citep[see][]{BerryHarlim2016, Bakry2013}. As a remedy, one performs various modifications to the kernel to adapt to these effects. We describe one such modification next.
	
	\paragraph{Bistochastic kernels} These are normalized / modified versions of $k_\epsilon$ which retain the symmetry along with additional properties such as Markov property. First, we define two functions
	\[ \deg_R(z) := \int k_\epsilon( z, z') d\mu(z'), \quad \deg_L(z') := \int \frac{ k_\epsilon( z, z') }{ \deg_R(z) } d\mu(z) . \]
	These are called the right and left degree functions respectively. Next, define a kernel
	\[ \tilde{k}_\epsilon(z, z'') := \frac{k_\epsilon(z, z'')}{ \deg_R(z) \deg_L(z'')^{1/2} } , \]
	and finally set
	\[ p_\epsilon(z, z') := \int \tilde{k}_\epsilon(z, z'') \tilde{k}_\epsilon(z'', z') d\mu(z'') . \]
	The kernel $p_\epsilon$ is symmetric and has stronger properties, revealed by considering its associated \emph{integral operator}. There are a number of variations of the bistochastic kernel \cite[e.g.][]{MarshallCoifman2019, WormellReich2021}, with different adaptations according to different contexts.
	
	\paragraph{Kernel integral operators} Given a $C^r$ kernel $\kappa : \Omega \times \Omega \to \real$, the associated integral operator $\mathcal{K}$ operates on $L^2(\mu)$ functions $\phi:\Omega\to \real$ as
	\[ (\mathcal{K}  \phi)(z) := \int_{\Omega} \kappa(z, z') \phi(z') d\mu(z'),  \]
	leading to a $C^r$ function $\mathcal{K}  \phi$. The integral operators corresponding to the kernels $\tilde{k}_\epsilon$ and $p_\epsilon$ are denoted as $\tilde{K}_\epsilon$ and $P_\epsilon$ respectively, defined similarly as 
	\[ (\tilde{K}_\epsilon  \phi)(z) := \int_{\Omega} \tilde{k}_\epsilon(z, z') \phi(z') d\mu(z'), \quad (P_\epsilon  \phi)(z) := \int_{\Omega} p_\epsilon(z, z') \phi(z') d\mu(z').  \]
	The kernels $\tilde{k}_\epsilon$ and $p_\epsilon$ have been designed so that we have the relation 
	\[ P_\epsilon = \tilde{K}_\epsilon \tilde{K}_\epsilon^* , \]
	where $\tilde{K}_\epsilon^*$ is the adjoint with respect to the Hilbert space structure of $L^2(\mu)$. It is a well known fact from Analysis that $ P_\epsilon$ is a  \emph{compact}, symmetric operator on $L^2(\mu)$ \citep[e,g,][]{DGJ_compactV_2018}. Moreover, $ P_\epsilon$ has a complete basis of unit-norm eigenfunctions
	\[  P_\epsilon \phi_j = \lambda_j \phi_j , \quad j= 1,2,\ldots , \]
	where the indexing is done so that the $\lambda_j$s are in decreasing order. Due to the bistochastic normalization, we have $\phi_1 \equiv 1_\Omega$, the constant function equal to $1$ everywhere. Moreover,  the eigenvalues satisfy $1 = \lambda_1 \geq \lambda_2 \geq \lambda_2 \geq \ldots > 0$. By self-adjointness of $P_\epsilon$, the $\phi_j$ form an orthonormal basis, i.e.,
	\[ \langle \phi_i, \phi_j \rangle_{L^2(\mu)} := \int \phi_i^*(x) \phi_j(x) d\mu(x) = \delta_{i,j} , \]
	where the $^*$ denotes a complex conjugate. All these properties of the $\lambda_j$ and $\phi_j$ are useful for \emph{kernel-based learning}, in which we recreate or extrapolate an unknown function from some samples, using these $\phi_j$s as a basis.
	
	\paragraph{Kernel based learning} Given an unknown function $f : \real^{\alpha(Q+1)} \to \real^d$, a \emph{learning} technique tries to approximate it by vectors chosen from some suitable function space, called a search / hypothesis space. In our case, this space is a finite-dimensional subspace of $L^2(\mu)$ spanned by eigenvectors $\phi_j$ of the kernel integral operator. A main advantage of a kernel based approach is that the $\phi_j$ can be approximated to any degree of accuracy by solving an eigenvalue equation of a data-driven matrix [see Algorithm~\ref{algo:kernel}]. A second advantage is that the $\phi_j$ extends easily from vectors to a continuous function over the entire data space $\real^{\alpha (Q+1)}$, irrespective of the nature of the measure $\mu$
	\begin{equation} \label{eqn:dskn30}
		\bar\phi_j(z) = \lambda_j^{-1}\int p_\epsilon(z , z') \phi_j(z') d\mu(z') = \lambda_j^{-1}\int \tilde{k}_\epsilon(z , z') \gamma_j(z') d\mu(z') .
	\end{equation}
	Note that the $\phi_j$ are members of the space $L^2(\mu)$. This makes the $\phi_j$ ill-defined as functions on the whole space $\real^{\alpha(Q+1)}$, or well-defined on only a subset of this space. Equation \eqref{eqn:dskn30} shows that due to the smoothing action of a kernel integral operator, $\phi_j$ extends to an actual function $\bar\phi_j$. It is an extension in the sense that 
	\[\bar\phi_j(z) = \phi_j(z), \quad \mu-a.e. \, z . \]
	We use a similar extension later in \eqref{eqn:def:barphi}, with $\mu$ replaced by a sampling measure. 
	
	In a data-driven setting, these integrals are replaced by matrix multiplications, see Algorithm~\ref{algo:oss} in Section~\ref{sec:data} for a precise description. 
	First, the components of $f$ along the $\phi_j$ are computed :
	\[ f_l := \langle \phi_l, f \rangle_{L^2(\mu)} = \int_{\real^{\alpha(Q+1)}} \phi_l^*(z) f(z) d\mu(z) ,  \]
	and then reconstructing $f$ as $f = \sum_l f_l \phi_l$. Since the $\phi_j$ form an orthonormal basis for $L^2(\mu)$, this reconstruction is possible in an $L^2(\mu)$ metric, for any $f\in L^2(\mu)$. In particular, this would be possible for a continuous function $f$ on $\real^{\alpha(Q+1)}$. Then by \eqref{eqn:dskn30} we have for every $x\in \real^{\alpha(Q+1)}$,
	\[\begin{split}
		f(x) &= \sum_l f_l \phi_l(x)  = \frac{1}{ \deg_R(x) } \int k_\epsilon(x, y) \sum_l f_l \tilde{\gamma}_l(y) d\mu(y) = \frac{ \langle k_x , w \rangle_{L^2(\mu)} }{ \langle k_x , 1 \rangle_{L^2(\mu)} }, \quad w := \sum_l f_l \frac{1}{\lambda_l} \gamma_l ,
	\end{split}\]
	where $k_x$ is the function $k(x,\cdot)$, known as the \emph{kernel section} at $x$. The function $w$ is called the \emph{feature vector} corresponding to the function $f$. The correspondence between $f$ and $w$ is linear and via an operator $Feature : L^2(\mu) \to L^2(\mu)$. The function $w = \Feature(f)$ plays the role of a density function. Note that we have a \emph{Markov transition function}
	\[ k^{\text{Markov}}_\epsilon : \Omega \times \Omega \to \real, \quad k^{\text{Markov}}_\epsilon (x,y) := k_\epsilon(x, y) / \int k_\epsilon (x, y') d\mu(y') . \]
	By design, for every $x\in \Omega$, $k^{\text{Markov}}_\epsilon(x,\cdot)$ is a probability density function on $\Omega$ with respect to $\mu$. The value of $f(x)$ then becomes the Markov transform of the initial distribution $w$ :
	\[ f(x) := \int k^{\text{Markov}}_\epsilon(x, y) w(y) d\mu(y). \]
	In practice, instead of doing an infinite sum of the form $\sum_l$, we use
	\[  f(y) \approx \sum_{l=1}^{L} f_l \phi_l(y) = \int_{\real^{\alpha(Q+1)}} p_\epsilon ( y, z) \sum_{l=1}^{L} \frac{ f_l }{ \lambda_l } \phi_l(z) d\mu(z) .  \]
	The parameter $L$ is called the \emph{spectral truncation} parameter. The higher the value of $L$, the more accurate the approximation is.  However higher order eigenfunctions of a matrix are more expensive to compute,  and their convergence to the limiting vector (in the sense of \cite{VonLuxburgEtAl2008}) is slower. The division by $\lambda_l$, which goes to zero as $l\to \infty$ also restricts how big $L$ could be, for a given computational resources and data-size.
	
	Note that our learning approach has no explicit null hypothesis on the map $g_{chaos}$, beyond the assumption that it is $C^2$ and thus square integrable. Under certain conditions, kernel eigenfunctions approximate Laplacian eigenfunctions \citep[e.g.][]{CoifmanLafon06b, BDGV_spectral_2020, VaughnBerryAntil2019} and tend to have zero derivatives when away from the dataset. This creates no conflict, as the kernel eigenfunctions still form a complete orthonormal basis on the support of the measure $\mu$. On the contrary, the vanishing of the derivatives away from the dataset guarantees that the simulated dynamics continues to have bounded orbits in the higher dimensional space $\real^{\alpha(Q+1)}$. 
	
	This completes a description of the theoretical basis of our methods. We next discuss the procedure in the data-driven setting of Assumption~\ref{A:data}. There, we also discuss the procedure for frequency identification in Algorithm~\ref{algo:RKHS_kernel}.
	
	\section{The data-driven procedure} \label{sec:data}
	
	In a data-driven approach, all of the spaces, operators and maps described in Section~\ref{sec:kernel} are approximated via data. This begins with an approximation of the dynamics-invariant measure $\mu$ by the \emph{sampling measure}
	\[\mu_N = \frac{1}{N} \sum_{n=1}^{N} \delta_{y_n^{(Q)}} , \]
	the average of the Dirac-delta measure on the data points $y_n^{(Q)}$. These approximate $\mu$ in a weak-sense via their integrals for every continuous test function $\phi : \real^{\alpha(Q+1)} \to \real$ :
	\[  \int_{\real^{\alpha(Q+1)}} \phi d \mu_N = \frac{1}{N} \sum_{n=1}^{N} \phi\left(  y_n^{(Q)} \right) \xrightarrow{n\to \infty} \int_{\real^{\alpha(Q+1)}} \phi d \mu . \]
	As a result, the infinite dimensional Hilbert space $L^2(\mu)$ will be represented as $L^2(\mu_N)$. The kernel integral operators $K$, $\tilde{K}$ and $P$ will be represented by $N\times N$ matrix $\Matrix{K}$, $\Matrix{\tilde K}$ and $\Matrix{P}$, as described below : 
	
	\begin{algorithm}[Kernel building] \label{algo:kernel}
		\ 
		\begin{enumerate}
			\item Input : 
			\begin{enumerate} [(i)]
				\item Data $\SetDef{ y_n\in \real^{\alpha}}{ n=1,\ldots,N }$ as in Assumption~\ref{A:data}.
				\item Bandwidth parameter $\epsilon>0$ for the kernel.
				\item The number of eigenvectors $L\in\num$ to be computed.
				\item Number of delay-coordinates $Q$.
			\end{enumerate}
			\item Output
			\begin{enumerate} [(i)]
				\item Eigenfunctions $\SetDef{ \vec{\phi}_{l} }{ l=1,\ldots, L }$ and eigenvalues $1=\lambda_1 \geq \ldots \geq \lambda_L$ of a bistochastic, symmetric kernel $p_\epsilon(x,y)$.
				\item Right singular vectors $\SetDef{ \vec{\gamma}_{l} }{ l=1,\ldots, L }$.
				\item Right degree vector $\vec{d}_r \in \real^N$.
			\end{enumerate}
			\item Steps
			\begin{enumerate} [(i)]
				\item Compute a $N\times N$ kernel matrix $\Matrix{K}$ using the Gaussian kernel $k_\epsilon$ \eqref{eqn:def:Gauss_ker} as :
				\[  \Matrix{K}_{i,j} := k_\epsilon \left(  y_i^{(Q)}, y_j^{(Q)} \right), \quad 1\leq i,j \leq N. \]
				\item Compute the degree vectors
				\[  \vec{d}_l := \frac{1}{N}  \Matrix{K} \vec{1}_N, \quad D_l:= \diag(  \vec{d}_l  ), \quad  \vec{d}_r := \frac{1}{N} \Matrix{K} D_l^{-1}, \]
				and then the matrix 
				\[ \Matrix{\tilde{K}} := D_l^{-1} \Matrix{K} D_r^{-0.5}, \quad D_r := \diag( \vec{d}_r ) . \]
				\item Compute the top $L$ singular values $1=\sigma_1 \geq \ldots \geq \sigma_L$ of $\Matrix{\tilde{K}}$ and the corresponding left eigenvectors $ \vec{\phi}_{1}, \ldots,  \vec{\phi}_{L}$ and right singular vectors $\vec{\gamma}_{1}, \ldots, \vec{\gamma}_{L}$.
				\item Set $\lambda_i := \sigma_i^2$, for $i=1,\ldots, L$.
			\end{enumerate}
		\end{enumerate}
	\end{algorithm}
	
	Algorithm~\ref{algo:kernel} is an initial processing step on the data. It is not specific to the reconstruction problem for the dynamics. The role of the bandwidth $\epsilon$ in numerical experiments is discussed in Section~\ref{sec:examples}. The output of the algorithm can be used for any learning problem based on the given dataset. The set of vectors $\SetDef{ \vec{\phi}_{l} }{ l = 1, \ldots, L }$ and $\SetDef{ \vec{\gamma}_{l} }{ l = 1, \ldots, L }$ are both orthonormal systems for $\cmplx^N$. The $N$-dimensional vectors $\vec{\phi}_{l}$ have continuous extensions to the whole of $\real^{\alpha(Q+1)}$ as
	\begin{equation} \label{eqn:def:barphi}
		\bar{\phi}_l : \real^{\alpha(Q+1)} \to \real, \quad  \bar{\phi}_l (y) := \frac{1}{N \lambda_l} \vec{ \tilde{k} }_\epsilon (y) ^\top \vec\gamma_l , \quad \vec{ \tilde{k} }_\epsilon (y) := \left(   \tilde{k}_\epsilon \left( y, y_1^{(Q)} \right)  , \ldots, \tilde{k}_\epsilon \left( y, y_N^{(Q)} \right) \right).
	\end{equation}
	The function $ \bar{\phi}_l$ is as smooth as the kernel $k_\epsilon$. If $y$ in the above equation is substituted by one of the data-points $y_n^{(Q)}$, then by design 
	\[ \bar{\phi}_l \left( y_n^{(Q)} \right) = \vec{\phi}_{l,n} . \]
	Thus $ \bar{\phi}_l$ is indeed a continuous extension of the vector $\vec{\phi}_l$. This feature of extendability and easy evaluation at arbitrary points is one of the most powerful tools of kernel-based methods. In our next algorithm, we show a different application of these eigenfunctions, for discovering true Koopman eigenfrequencies. It is based on ergodic theoretic results derived in \cite{DasGiannakis_RKHS_2018}. It involves using the familiar fast-Fourier transform but on the $L$ eigenfunctions derived above instead of the raw data, along with a weighting using the $\lambda_l$s. The results are interpreted in a functional space called a \emph{reproducing kernel Hilbert space} or RKHS.
	
	\begin{algorithm}[RKHS based spectral analysis] \label{algo:RKHS_kernel}
		\citep[][Alg. 1]{DasGiannakis_RKHS_2018}.
		\begin{enumerate}
			\item Input 
			\begin{enumerate} [(i)]
				\item Eigenfunctions $\SetDef{ \vec{\phi}_{l} }{ l=1,\ldots, L }$ and eigenvalues $1=\lambda_1 \geq \ldots \geq \lambda_L$ from Algorithm~\ref{algo:kernel}.
				\item Threshold parameters : $\epsilon_1, \epsilon_2>0$ and integer $L_0$ such that $1<L_0<L$.
				\item Sampling interval $\Delta t$ if the source is a continuous time system.
			\end{enumerate}
			\item Output : A set of frequencies identified $0=\omega_1<\omega_2<\ldots<\omega_m$ identified as true Koopman eigenfrequencies.
			\item Steps
			\begin{enumerate} [(i)]
				\item Collect the eigenvectors $\vec{\phi}_{l}$ in an $N\times L$ matrix $\Matrix{\Phi}$. Let $\Fourier_N$ be the discrete Fourier transform on $N$ vectos. Set $\Lambda:= \diag\left(  \lambda_1, \ldots,\lambda_L \right)$ and compute
				\[ \Matrix{\hat\Phi} := \Fourier_N \Matrix{\Phi} \quad \Matrix{H} := \Matrix{\hat\Phi} \Lambda^{-0.5}. \]
				\item Next compute an $N\times L$ matrix $\Matrix{W}$ such that for each $n=1,\ldots, N$,
				\[ \Matrix{W}_{n,1} := \abs{ \Matrix{H}_{n,1} }, \quad  \Matrix{W}_{n,l+1} := \Matrix{W}_{n,l} + \abs{ \Matrix{H}_{n,l+1} }, \quad l=1,\ldots, L-1 .  \]
				\item Set $J = 1,\ldots, N$. 
				\item Discard all the $j\in J$ for which $\Matrix{W}_{j,L_0}<\epsilon_1$.
				\item Of the remaining $j\in J$, discard those $j$ for which $\ln \Matrix{W}_{j,L} - \ln \Matrix{W}_{j,L_0} > \epsilon_2$.
				\item Compute $\omega_j = \frac{2\pi j}{N}$ for all of the remaining $j\in J$.
				\item If the underlying system is continuous time, then divide each of the $\omega_j$ by $\Delta t$.
			\end{enumerate}
		\end{enumerate}
	\end{algorithm}
	
	Note that there are two filterings taking place, in steps (iii) and (iv), via parameters $\epsilon_1$ and $\epsilon_2$ respectively. They are based on results in approximation theory on Reproducing kernel Hilbert spaces \citep[see][Thm 1, 4]{DasGiannakis_RKHS_2018}. The identified frequencies $0=\omega_1<\omega_2<\ldots<\omega_m$ are by no means exhaustive, they are only a finite subset of a usually infinite set of Koopman eigenfrequencies. However, they represent those (true) frequencies that have a significant presence in the data. The threshold $\epsilon_1$ is meant to be a numerical implementation of frequencies being significant. We shall use these selected frequencies later to build our reconstructed dynamics \eqref{algo:oss}. Algorithm~\ref{algo:RKHS_kernel} is unique in its use of RKHS-regularity as a criterion for identifying frequencies. See Table~\ref{tab:compare} for a comparison of our methods with other techniques. 
	
	We utilize the delay-coordinate structure in the embedding to simplify the construction of $g_{chaos} : \real^{\alpha(Q+1)} \to \real^{\alpha}$ as
	\begin{equation} \label{eqn:def:gchaos}
		g_{chaos} \left( x^{(0)} , \ldots, x^{(Q)} \right) := \Matrix{\begin{array}{c}
				\hat{g}_{chaos}^{(0)} \left( x^{(0)} , \ldots, x^{(Q)} \right)  \\
				x^{(0)} \\
				\vdots \\
				x^{(Q-1)}
		\end{array}}
	\end{equation}
	Here for each $q\in 0,\ldots, Q$, $x^{(q)}$is the $q$-th set of coordinates in $\real^{\alpha}$. We next describe reconstruct the periodic and chaotic components $g_{per}$ and $\hat{g}_{chaos}^{(0)}$. The set of Koopman eigenfrequencies  identified from Algorithm~\ref{algo:RKHS_kernel} can be passed as the second input to the Algorithm below. Given a matrix $\Matrix{E}$, we use $\Matrix{E}_{l,:}$ to denote its $l$-th row. Also recall the functions $\bar{\phi}_l$ defined in \eqref{eqn:def:barphi}.
	
	\begin{algorithm}[Components of the dynamics] \label{algo:ker_chaos}
		\  
		\begin{enumerate}
			\item Input
			\begin{enumerate} [(i)]
				\item An $N\times \alpha$ matrix $\Matrix{Y}$ whoose $n$-th row represents the data point $y_n\in \real^{\alpha}$ as in Assumption~\ref{A:data}. 
				\item A set of Koopman eigenfrequencies  $0=\omega_1<\omega_2<\ldots<\omega_m$.
				\item $N\times L$ matrix $\Phi$ whose columns store the eigenfunctions $\vec{\phi}_l$ from Algorithm~\ref{algo:kernel}.
			\end{enumerate}
			\item Output : 
			\begin{enumerate} [(i)] 
				\item $L\times \alpha$ matrix $\Matrix{E}$ which approximates the function
				\begin{equation} \label{eqn:def:Ecoeff}
					\hat{g}_{chaos}^{(0)} : \real^{\alpha(Q+1)} \to \real^{\alpha}, \quad \hat{g}_{chaos}^{(0)} \approx \sum_{l=1}^{L} \Matrix{E}_{l,:} \bar{\phi}_l .
				\end{equation}
				\item $m\times \alpha$ matrix $\Matrix{A}$ which represents a periodic function 
				\begin{equation} \label{eqn:def:Acoeff}
					\hat{g}_{per} : \real \to \real^{\alpha}, \quad g_{per}(t) := \Re \sum_{j=1}^{m} \Matrix{A}_{j,:} e^{\iota \omega_j t} .
				\end{equation}
			\end{enumerate}
			\item Steps
			\begin{enumerate} [(i)]
				\item Define an $N\times m$ matrix $\Matrix{F}$ as
				\[ \Matrix{F}_{n,j} := (2 - \delta_{j,1}) e^{\iota n \omega_j}  \quad 1\leq n\leq N, \, 1\leq j\leq m .  \]
				\item Find an $m\times \alpha$ matrix $A$ which is the least-squares solution to 
				\[ \Re \Matrix{F} \Matrix{A} = \Matrix{Y} \]
				\item Set $\Matrix{Y_{non}} :=  \Matrix{Y} - \Re \Matrix{F} \Matrix{A}$ and the $L\times \alpha$ matrix $\Matrix{E} := \Matrix{\Phi}^* \Matrix{Y_{non}}$.
			\end{enumerate}
		\end{enumerate}
	\end{algorithm}
	
	Algorithms \ref{algo:kernel}, \ref{algo:RKHS_kernel} and \ref{algo:ker_chaos} reconstruct the quasiperiodically driven dynamics which underlies the data. We shall describe another algorithm to perform the evaluations of the functions $\bar{\phi}_l$ as in \eqref{eqn:def:Ecoeff}.
	
	\begin{algorithm}[Out of sample evaluations] \label{algo:oss}
		\ 
		\begin{enumerate}
			\item Input
			\begin{enumerate} [(i)]
				\item Data $\SetDef{ y_n\in \real^{\alpha}}{ n=1,\ldots,N }$ as in Assumption~\ref{A:data}, along with sampling interval $\Delta t$ if the source is a continuous time system.
				\item $L\times \alpha$ matrix $\Matrix{E}$ from Algorithm \ref{algo:ker_chaos}.
				\item Vector $y\in \real^{\alpha(Q+1)}$ representing a point of evaluation.
				\item Right degree vector $\vec{d}_r$ from Algorithm \ref{algo:kernel}.
				\item Right singular vectors $\SetDef{ \vec{\gamma}_{l} }{ l=1,\ldots, L }$ from Algorithm \ref{algo:kernel}.
			\end{enumerate}
			\item Output : vector $\hat{g}_{chaos}^{(0)}(y) \in \real^{\alpha}$. 
			\item Steps
			\begin{enumerate} [(i)]
				\item Create an $N\times L$ matrix $\Matrix{ \tilde{\Gamma} }$ defined as $\Matrix{ \tilde{\Gamma} }_{n,l} := \paran{ \vec\gamma_l }_{n} \lambda_l^{-1/2} \paran{ \vec{d}_r }_n^{-1}$.
				\item Next compute
				\[ \vec{k_{os}} := \left(  k_\epsilon \left( y_n^{(Q)}, y  \right) \right)_{n=1}^{N}, \quad s :=  \frac{1}{N} \sum_{n=1}^{N}  \vec{k_{os}}(n).  \]
				\item Finally compute 
				\[  \left( e_1(y), \ldots, e_k(y) \right) \approx \frac{1}{Ns} \vec{k_{os}}^\top \Matrix{ \tilde{\Gamma} } \Matrix{E} . \]
			\end{enumerate}
		\end{enumerate}
	\end{algorithm}
	
	\paragraph{The reconstruction} The $m$ eigenfrequencies identified by Algorithm~\ref{algo:RKHS_kernel} are generated by some collection of $d$ eigenfrequencies $\vec{ \tilde{\omega} } = \paran{ \tilde{\omega}_1, \ldots, \tilde{\omega}_d }$. Each selected frequency $\omega$ is thus of the form $\omega = \vec{j} \cdot \vec{ \tilde{\omega} }$ for some $d$-dimensional integer vector $\vec{j} = \left( j_1, \ldots, j_d \right)$. Then the periodic part $g_{per}$ can be written as the Fourier series
	\[ g_{per}(\vec\theta) = \sum_{\vec{j} \in \integer^d} a_{\vec{j}} \exp\left( \iota  (\vec{j} \star \vec{ \tilde{\omega} }  ) \cdot \vec{\theta} \right) ,  \]
	where $\vec{j} \star \vec{ \tilde{\omega} }$ is the $d$-dimensional vector $\left(  j_1 \tilde{\omega} _1, \ldots, j_d \tilde{\omega} _d\right)$. The output of Algorithm~\ref{algo:RKHS_kernel} provides a means of avoiding the task of identifying $\vec{ \tilde{\omega} } $. The $\omega$ that are selected correspond to those indices $\vec{j}$ for which $\abs{ a_{\vec{j}} }$ is substantial. Thus we have the approximation
	\begin{equation} \label{eqn:sdkjn39}
		g_{per} \left( \vec\theta + n\vec{ \tilde{\omega} } \right) \approx \hat{g}_{per} \left( n \Delta t \right) := \sum_{j=1}^{m} a_{j} \exp\left( \iota \Delta t n \omega_j  \right) ., \quad n=0,1,2,\ldots . 
	\end{equation}
	Using this simplification~\ref{eqn:sdkjn39}, and the formulas in \eqref{eqn:def:Ecoeff} and \eqref{eqn:def:Acoeff}, we create the following data-driven model of the dynamics :

	\begin{equation} \label{eqn:reconstruct} 
		\boxed{ \begin{split}
				y^0_{n+1} &:=  \hat{g}_{per}(n) + \hat{g}_{chaos}^{(0)}\left( y^0_n , \ldots, y^Q_n   \right) \\
				y^1_{n+1} &:= y^0_n \\
				\vdots &= \vdots \\
				y^Q_{n+1} &:= y^{Q-1}_n 
		\end{split} }, \quad y^q_n \in \real^{\alpha}, \quad q\in 0,\ldots, Q, \quad n=1,2,3,\ldots .
	\end{equation}

	\begin{theorem} [Convergence] \label{thm:2}
		Suppose Assumption~\ref{A:basic}, \ref{A:cnstnt_snstv} and \ref{A:data} hold, and assume that the number of delays $Q$ is large enough so that $Y^{(Q)} : \Omega \to \real^{ \alpha (Q+1) }$ is an injective map. Fix  an error bound $\gamma>0$ and a forecast time $T\in \num$. Then for large enough number of data samples $N$, and parameters $\epsilon_1, \epsilon_2$ from Algorithm \ref{algo:RKHS_kernel} small enough, and $L_0,L$ large enough, we have for $\nu$-almost every $x_0$ in $X$ :
		\[ \norm{ Y\circ g_{chaos} -  \hat{g}_{chaos}^{(0)} \circ Y^{(Q)} }_{L^2(\mu)} \leq \gamma; \quad \abs{ g_{per} \paran{ \theta_n} - \hat{g}_{per}(n\Delta t) } \leq \gamma, \quad \forall 0 \leq n \leq T  . \]
	\end{theorem}
	
	\begin{proof} Since $g$ is $C^2$, $g_{per}$ is $C^2$ too. As a result, the Fourier series of $g_{per}$ :
		\[ g_{per} (\theta) = \sum_{ \vec{j} \in \integer^d } a_{\vec{j}} \exp \paran{ \iota \vec{j} \cdot \theta } ,\]
		converges uniformly. The frequencies $\SetDef{ \vec{j} \cdot \vec{\omega } }{ \vec{j} \in \integer^d }$ can be ordered based on the lexical ordering of $\integer^d$. Thus there is an integer $m>0$ such if $\vec{j}^{(1)} , \ldots , \vec{j}^{(m)}$ are the first $m$ eigenfrequencies, then
		\[ \sup_{\theta \in \mathbb{T}^d} \abs{ g_{per} (\theta) - \sum_{ i=1 }^{m} a_{\vec{j}^{(i)}} \exp \paran{ \iota \vec{j}^{(i)} \cdot \theta } } \leq \gamma . \]
		The corresponding selected frequencies are $\omega_i = \vec{j}^{(i)} \cdot \rho$, for $i\in 1, \ldots, m$. According to \citep[][Corr 2]{DasGiannakis_RKHS_2018}, if the filtration parameters $\epsilon_1, \epsilon_2$ are chosen small enough, and number of kernel eigenfunction  $L_0, L$ chosen large enough, then there is an $N_1\in\num$ such that if the number of data-samples $N$ is greater than $N_1$, the first $m$ eigenfunctions $\zeta_1, \ldots, \zeta_m$ can be calculated within a precision of $\frac{2\pi}{N \Delta t}$, where $\Delta t$ is the sampling interval. As a result, the non-quasiperiodic component $Y_{non}$ in Algorithm \ref{algo:ker_chaos} is computed with an error that is less than $\gamma$. By the universal approximation property of RKHSs \citep[e.g.][]{Sriperumbudur_2009_kernel, HeinEtAl2005, VonLuxburgEtAl2008, ScholkopfSmolaMu1998}, $g_{chaos}$  can be approximated within an $L^2$ error of less than $\gamma$, when $N>N_2$ for some lower bound $N_2$. Thus any $N$ greater than $\max( N_1, N_2)$ achieves the desired approximation.
	\end{proof}
	
	This completes the statement of our methods and theory. We next show that the hypothesis of quasiperiodically driven dynamics proves effective for analyzing data from several real-world systems.
	
	\section{Case studies} \label{sec:examples}
	
	We analyze data from three real-world systems using our techniques from the previous sections.
	
	\begin{table}
		\caption{Summary of parameters in experiments. In all these experiments, $N=2\times 10^4$ and $\epsilon_1=0.1$. }
		\begin{tabularx}{\linewidth}{|L|L|L|L|L|L|L|L|}
			\hline
			Data-source & Figures & $\epsilon$ & $Q$ & $L_0$ & $\epsilon_2$ \\ \hline
			Freeway & \ref{fig:408_weekly}, \ref{fig:408_daily}, \ref{fig:408_eigen} & 0.3 & 50 & 100 & 3.1 \\ \hline
			Signalized intersections &\ref{fig:Alafaya_reconstruct}, \ref{fig:Alafaya_daily}, \ref{fig:Alafaya_eigen}, & 0.1 & 50 & 100 & 4.0 \\ \hline
			Heart atrial data & \ref{fig:atrial} & 0.002 & 50 & 500 & 5.1 \\ \hline
		\end{tabularx}
		\label{tab:param}
	\end{table}
	
\begin{figure}\center
		\includegraphics[width=.9\linewidth]{\figs 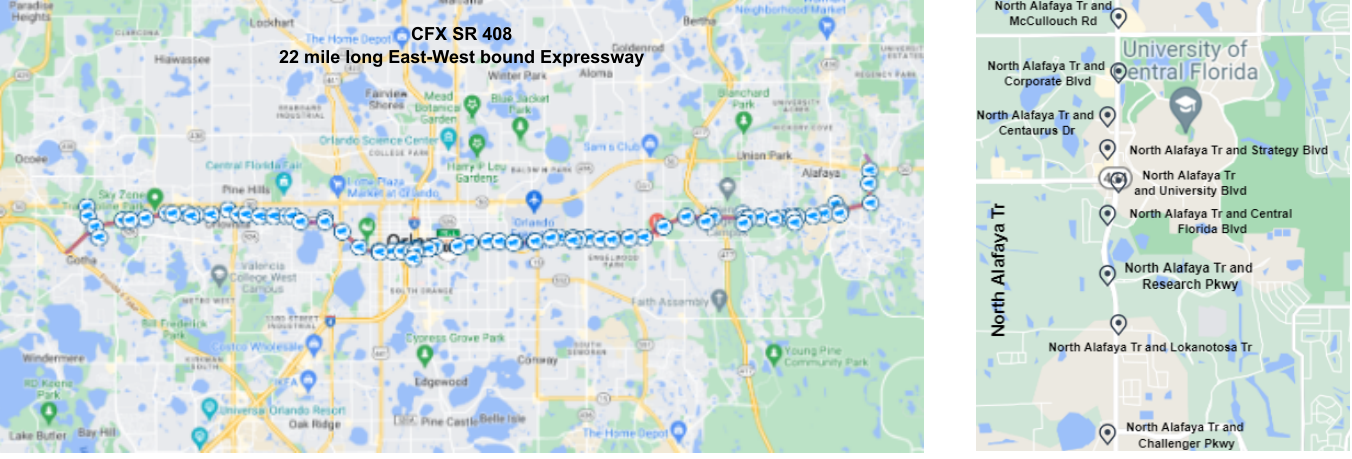}
		\caption{Location of CFX SR 408 freeway (left) and intersections on Alafaya corridor (right). Locations of cameras placed on SR 408 are shown on the map. The
			figures are taken from Google map}
		\label{fig:map}
	\end{figure}

\begin{figure}\center
		\includegraphics[width=.9\linewidth]{\figs 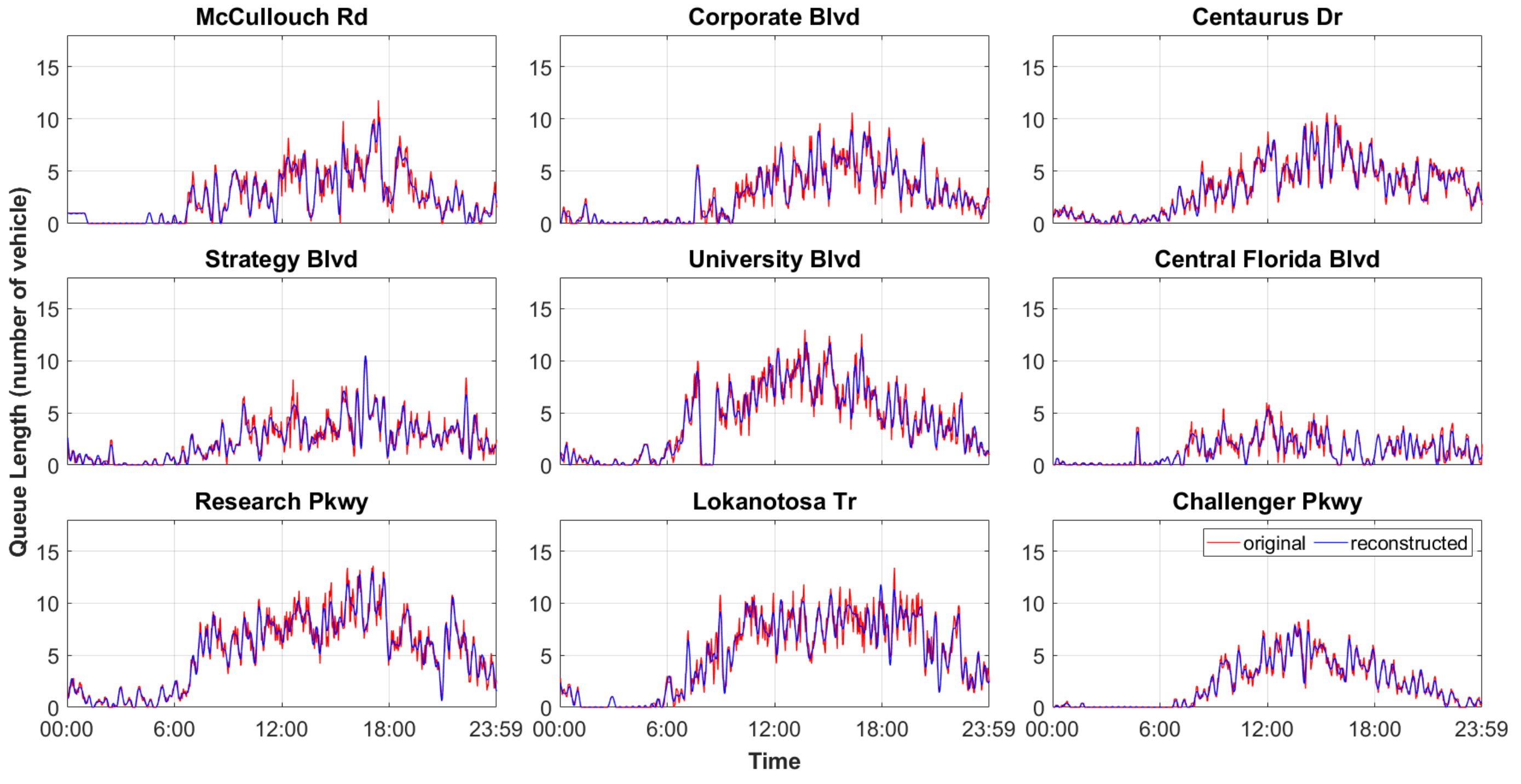}
		\caption{Closer look at the reconstruction of signalized intersection queue-length dynamics. The results shown in Figure~\ref{fig:Alafaya_reconstruct} are examined on a smaller time-scale of a day. The reconstruction appears accurate and coherent with the original traffic data.}
		\label{fig:Alafaya_daily}
	\end{figure}

\begin{figure}\center
		\includegraphics[width=.9\linewidth]{\figs 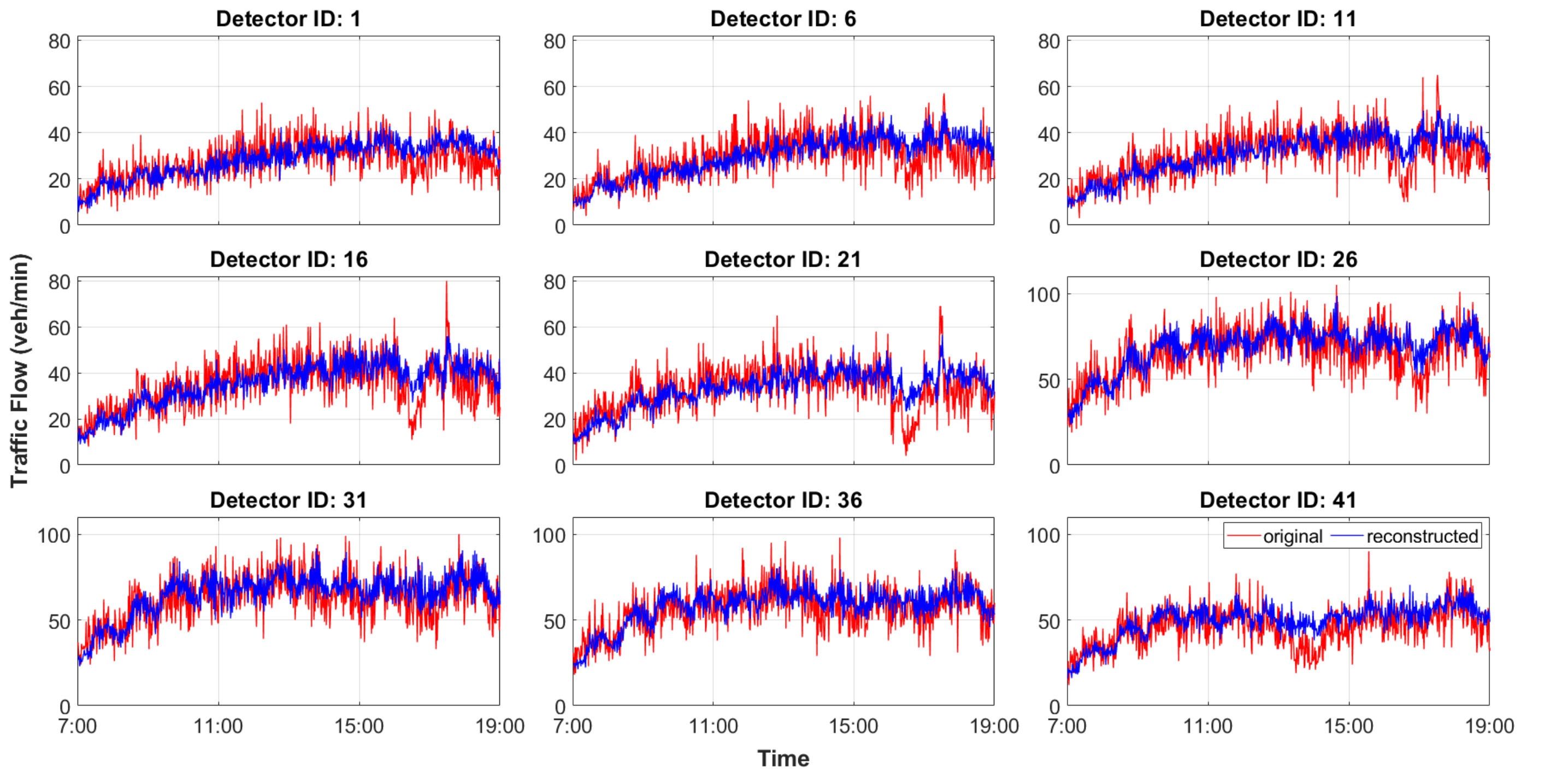}
		\caption{ Closer look at the reconstruction of freeway traffic flow dynamics. The results shown in Figure~\ref{fig:408_weekly} is examined on a smaller time-scale of a day. The reconstruction appears accurate and coherent with the original traffic data.}
		\label{fig:408_daily}
	\end{figure}

	\begin{enumerate}
		
		\item Traffic dynamics along signalized intersections. Queue length i.e. the number of stopped vehicles  at an intersection during red lights is an indicator of the dynamics in an urban area. In this study we obtained queue length data from $9$ signalized intersections of the Alafaya corridor situated in Orlando, Florida, United States. The location of the corridor on map is shown in Figure~\ref{fig:map} (right). This proprietary dataset was provided by \href{https://trafficbot.rhythmtraffic.com/}{InSync} and first used in \cite{rahman2021real}. The dataset contains the queue lengths aggregated over every 2 min which is the approximate traffic-signal cycle duration. 
		See Figure \ref{fig:Alafaya_reconstruct} for a view of the data, and its reconstruction, and Figure \ref{fig:Alafaya_daily} for a more magnified view.
		
		\item Traffic dynamics on the freeway. Traffic flow is a critical parameter to understand freeway traffic dynamics. Traffic flow refers to the number of vehicles passing through any cross-section of a road. We obtained freeway traffic flow data from the SR 408, which is a tolled expressway situated in Orlando, Florida, United States. The location of the freeway on the map is shown in Figure~\ref{fig:map} (left).  
		We used proprietary  traffic flow data collected from 64 sensors placed along a 22-mile corridor. We used the dataset provided by 
		\href{https://www.cfxway.com/for-travelers/expressways/408/}{Central Florida Expressway Authority (CFX)} in which outliers and faulty detectors were corrected.
		See Figure \ref{fig:408_weekly} for a view of the data, and its reconstruction, and Figure \ref{fig:408_daily} for a more magnified view.

		\item Cardiac signals. The data was obtained from \href{ https://physionet.org/content/afdb/1.0.0/ }{PhysioNet}'s repository of medical data. The data was collected in a study on the use of R-R intervals for detecting atrial fibrillation \cite{moody1983new}. The data is a timeseries covering a time duration of 2 hours, sampled at 250 samples per second. The timeseries contained two columns representing unedited recordings of atrial fibrillation and atrial flutter, with 12-bit resolution over a range of ±10 milivolts.
	\end{enumerate}
	
\begin{figure}[!htbp]\center
		\includegraphics[width=.80\linewidth]{\figs 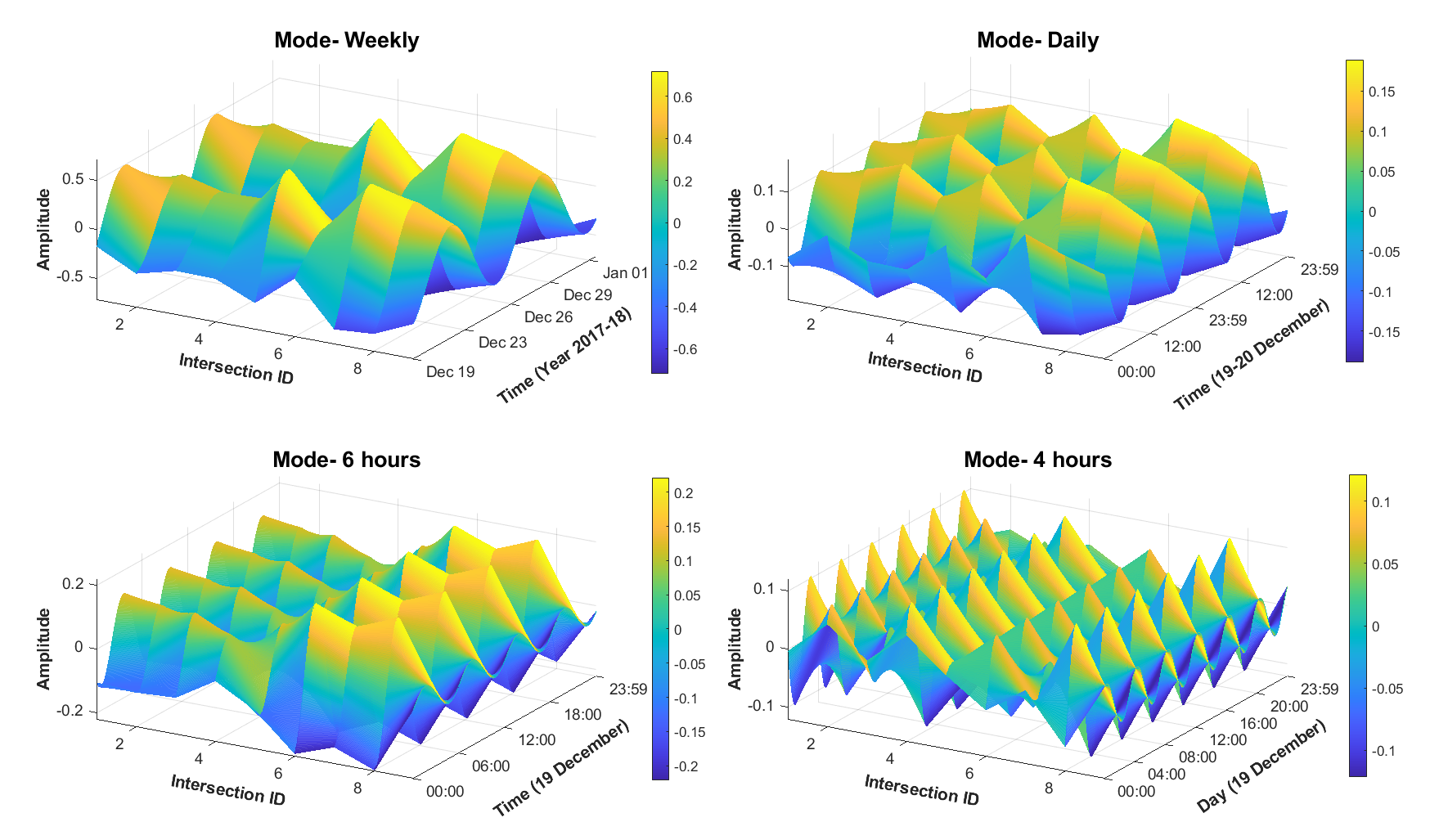}
		\caption{ Koopman eigenfunctions for signalized intersection dynamics. The traffic queue length build-up on the Alafaya corridor is assumed to be driven by a quasiperiodically driven dynamics, of the form \eqref{eqn:def:quasi_driven_II}. The Koopman eigenmodes of this dynamics were extracted using purely data-driven means (Algorithms \ref{algo:kernel}, \ref{algo:RKHS_kernel}). They describe coherent spatiotemporal patterns present within the dynamics. }
		\label{fig:Alafaya_eigen}
	\end{figure}

\begin{figure}\center
		\includegraphics[width=.98\linewidth]{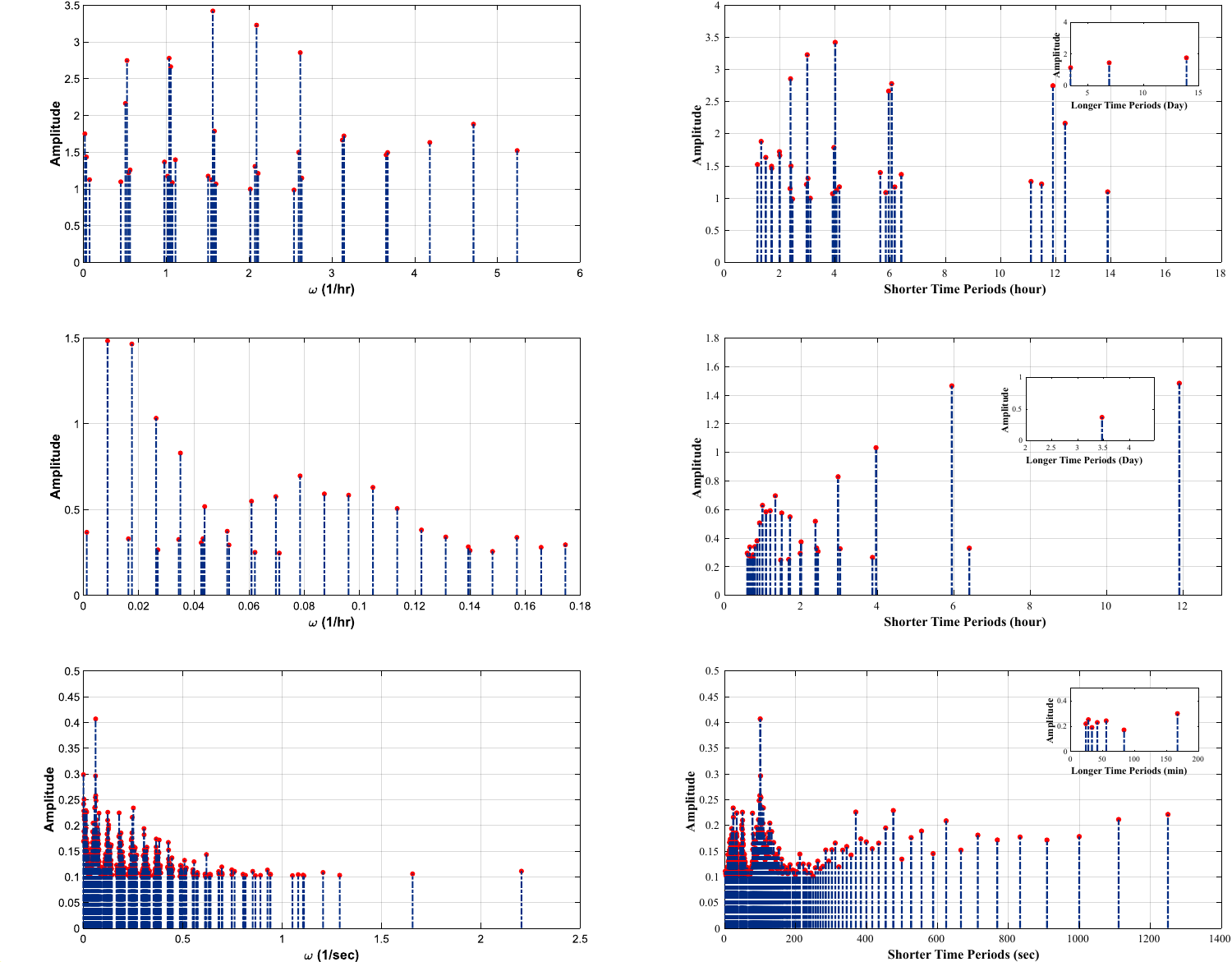}
		\caption{ Selected frequencies and time-periods. The top, middle, and bottom panels correspond to the signalized intersections (Alafaya corridor), freeway (CFX SR 408), and the cardiac data respectively. The identified frequencies correspond to Koopman frequencies of the dynamical system \eqref{eqn:def:quasi_driven_II} that we conceptually associate with the data. They can be interpreted as natural frequencies of these systems. Each time period $T$ that is marked corresponds to a selected frequency $\omega$ via the relation $T = 2\pi/ \omega \Delta T$, with $\Delta T$ being the sampling interval for the data. }
		\label{fig:freq_sel}
	\end{figure}

\begin{figure}\center
		
		\includegraphics[width=.32\linewidth]{\figs 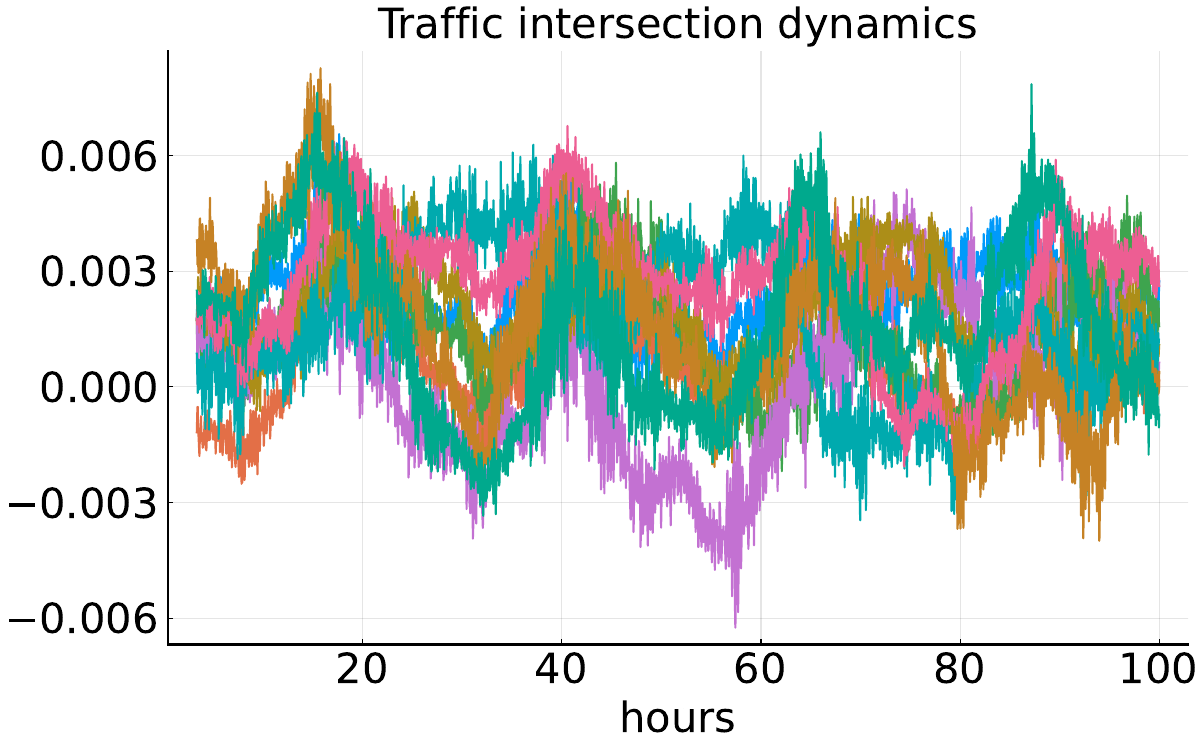}
		\includegraphics[width=.32\linewidth]{\figs 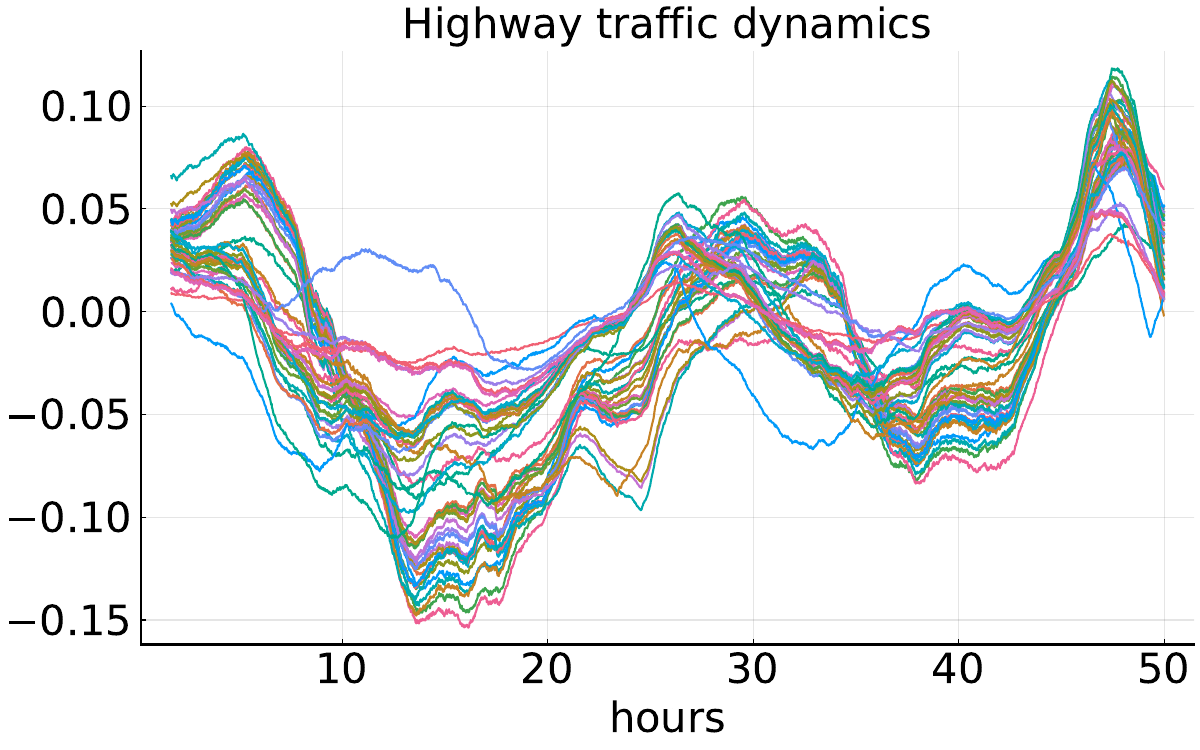}
		\includegraphics[width=.32\linewidth]{\figs 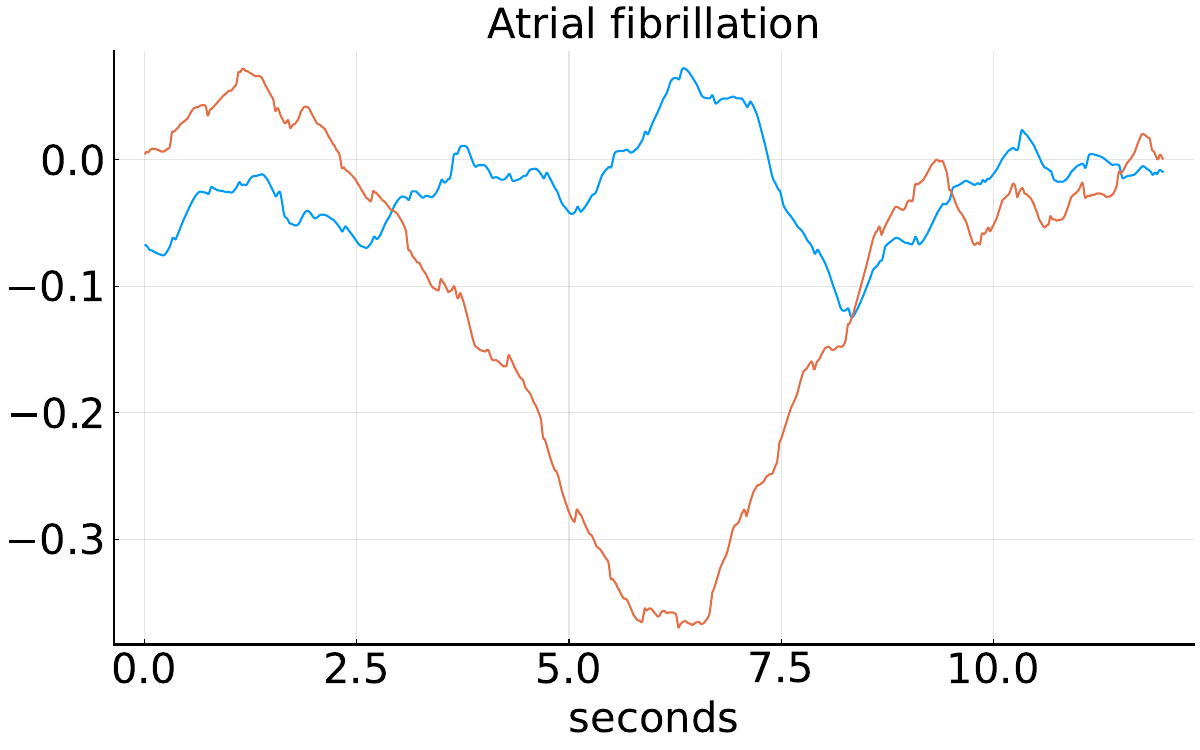}
		\caption{ Error analysis for the reconstructed dynamics. The amplitude normalized moving average error \eqref{eqn:def:anma_err} is computed for the three experiments, with a moving average window of $500$ time units. The number of curves in each plot corresponds to $\alpha$, the dimension of the signal used for analysis in that experiment. The use of a moving average window reduces the effect of random outliers or events. It diminishes the contribution of stochasticity to the error. Usually, when there is a mismatch between two dynamical systems such as the true dynamics \eqref{eqn:def:quasi_driven_II} and the reconstruction \eqref{eqn:reconstruct}, the two systems diverge and the error grows unrestricted. However, due to the special structure of quasiperiodically driven dynamics, and our ability to reconstruct the periodic part with high accuracy, our error does not grow monotonically but makes small oscillations.}
		\label{fig:error}
	\end{figure}

	Our numerical methods has three main objectives
	\begin{enumerate}[(i)]
		\item Reconstruction : Figures~\ref{fig:408_daily} and \ref{fig:Alafaya_daily} provide a zoomed view of the reconstruction plots. They reveal the same level of accuracy over smaller time-scales of a day instead of weeks. Given a reference timeseries $y_n$, and a reconstruction $\hat{y}_n$, the \emph{amplitude normalized, moving averaged error} is
		\begin{equation} \label{eqn:def:anma_err}
			e_{T,n} := \frac{1}{ \norm{y}_{\sup} } \frac{1}{T} \sum_{t=0}^{T-1} \SqBrack{ \hat{y}_{n+t} - y_{n+t}  }, \quad n=0,1,2,\ldots .
		\end{equation}
		See Figure \ref{fig:error} for the calculation of these errors, and the rationale behind this choice. The magnitude of this error is about 0.7\%, 15\%, and 3\%  respectively for the signalized intersections, freeway, and Cardiac data experiment.
		\item Identifying frequencies : Natural frequencies of the system contribute to the basic understanding of dynamical systems. They reveal the timescales as well as the quasiperiodicity dimension. These are displayed in Figure \ref{fig:freq_sel}.
		\item Identifying Koopman modes : In systems such as traffic systems, one has a direct interpretation of some of the variables as spatial coordinates, which vary with time. Spatiotemporal patterns such as those in Figures \ref{fig:Alafaya_eigen} represent components of the dynamics which is periodic in a combined space-time coordinate system.
	\end{enumerate}
	
	\begin{figure}\center
		\includegraphics[width=.32\linewidth]{\figs 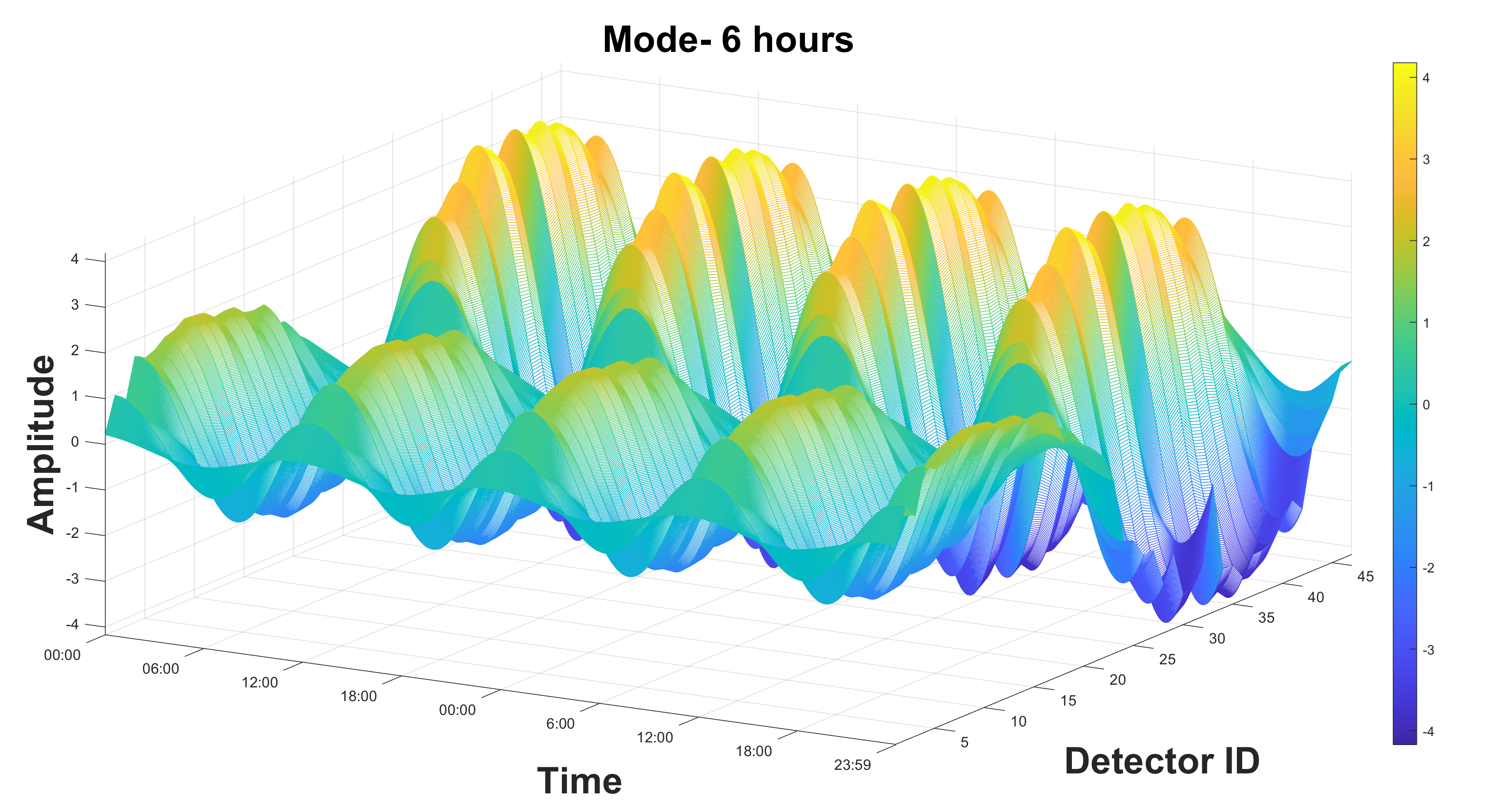}
		\includegraphics[width=.32\linewidth]{\figs 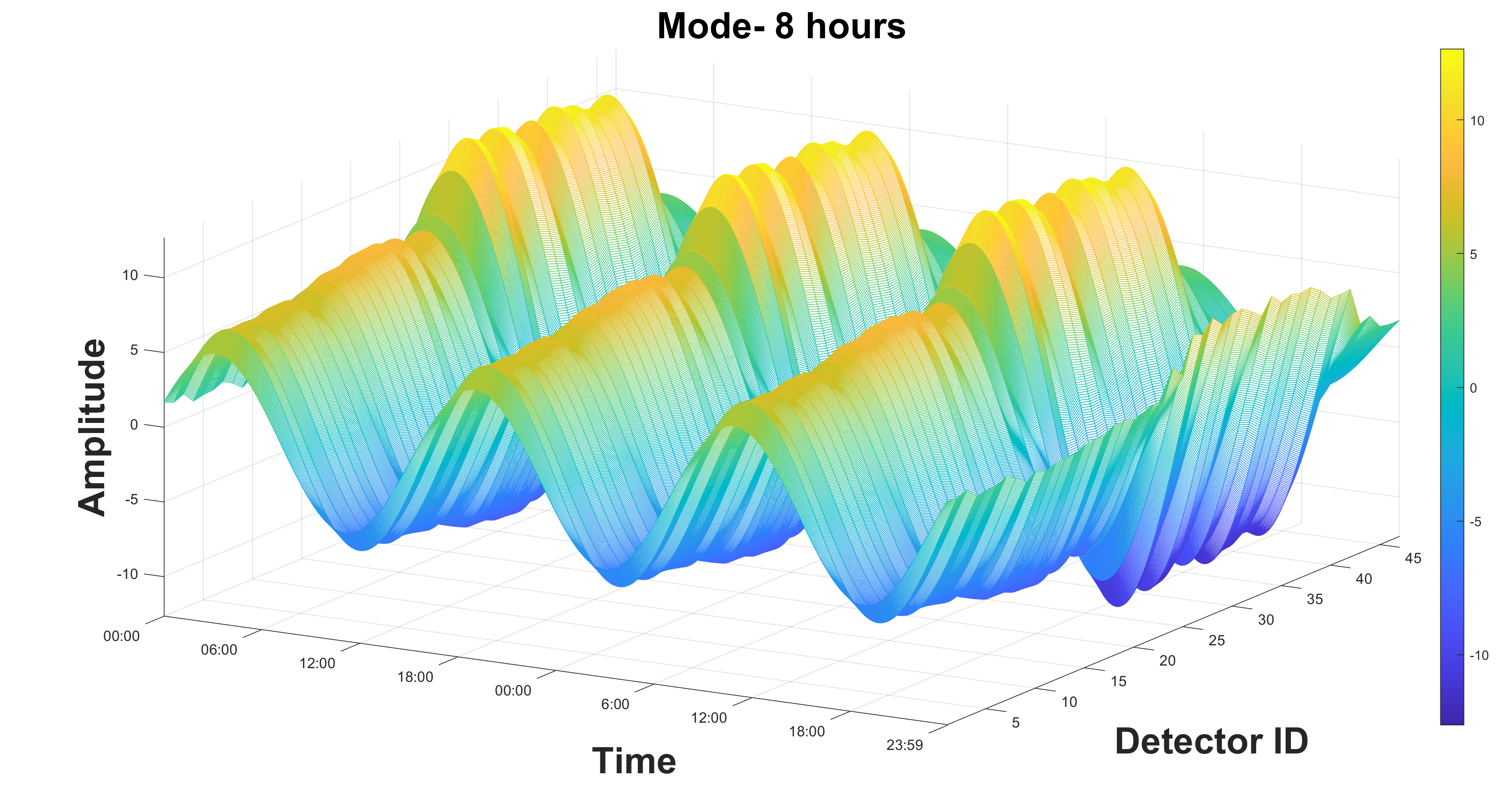}
		\includegraphics[width=.32\linewidth]{\figs 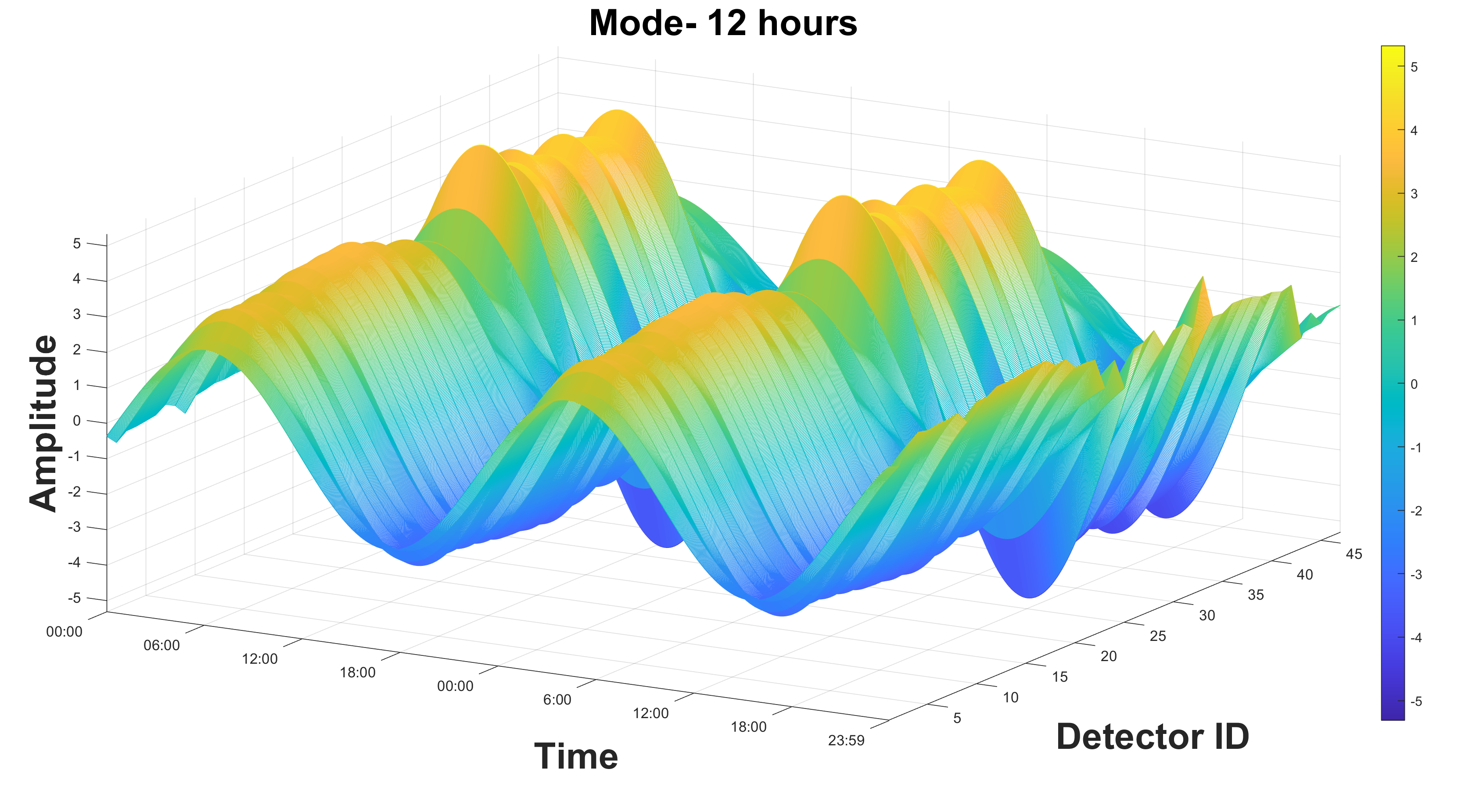}
		\caption{ Koopman eigenfunction for freeway traffic dynamics. The flux of traffic flow along the Florida 408 freeway is assumed to be driven by a quasiperiodically driven dynamics, of the form \eqref{eqn:def:quasi_driven_II}. The Koopman eigenmodes of this dynamics were extracted using purely data-driven means (Algorithms \ref{algo:kernel}, \ref{algo:RKHS_kernel}). They describe coherent spatiotemporal patterns present within the dynamics. }
		\label{fig:408_eigen}
	\end{figure}
	
	\begin{figure}\center
		\includegraphics[width=.31\linewidth]{\figs 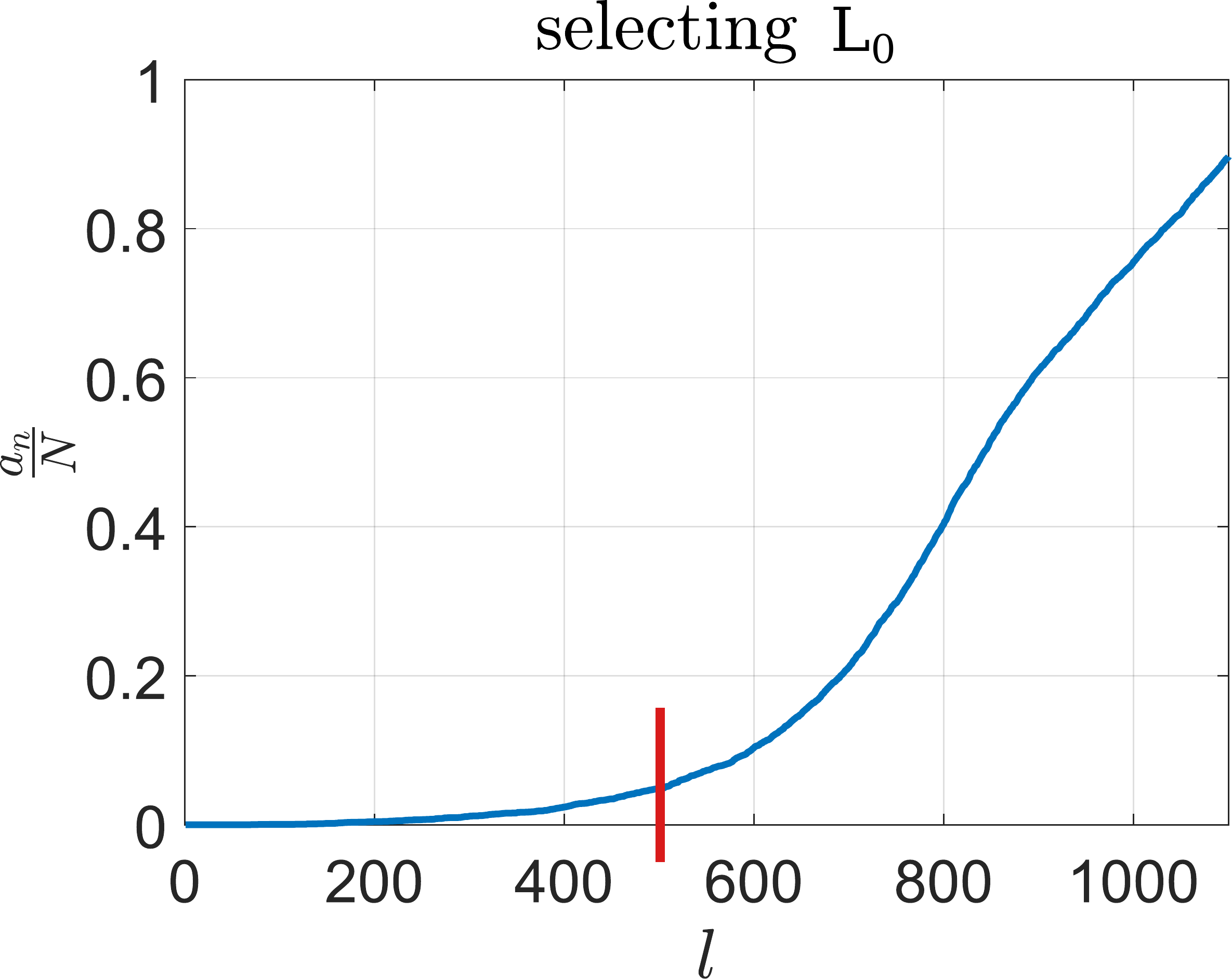}	
		\includegraphics[width=.31\linewidth]{\figs 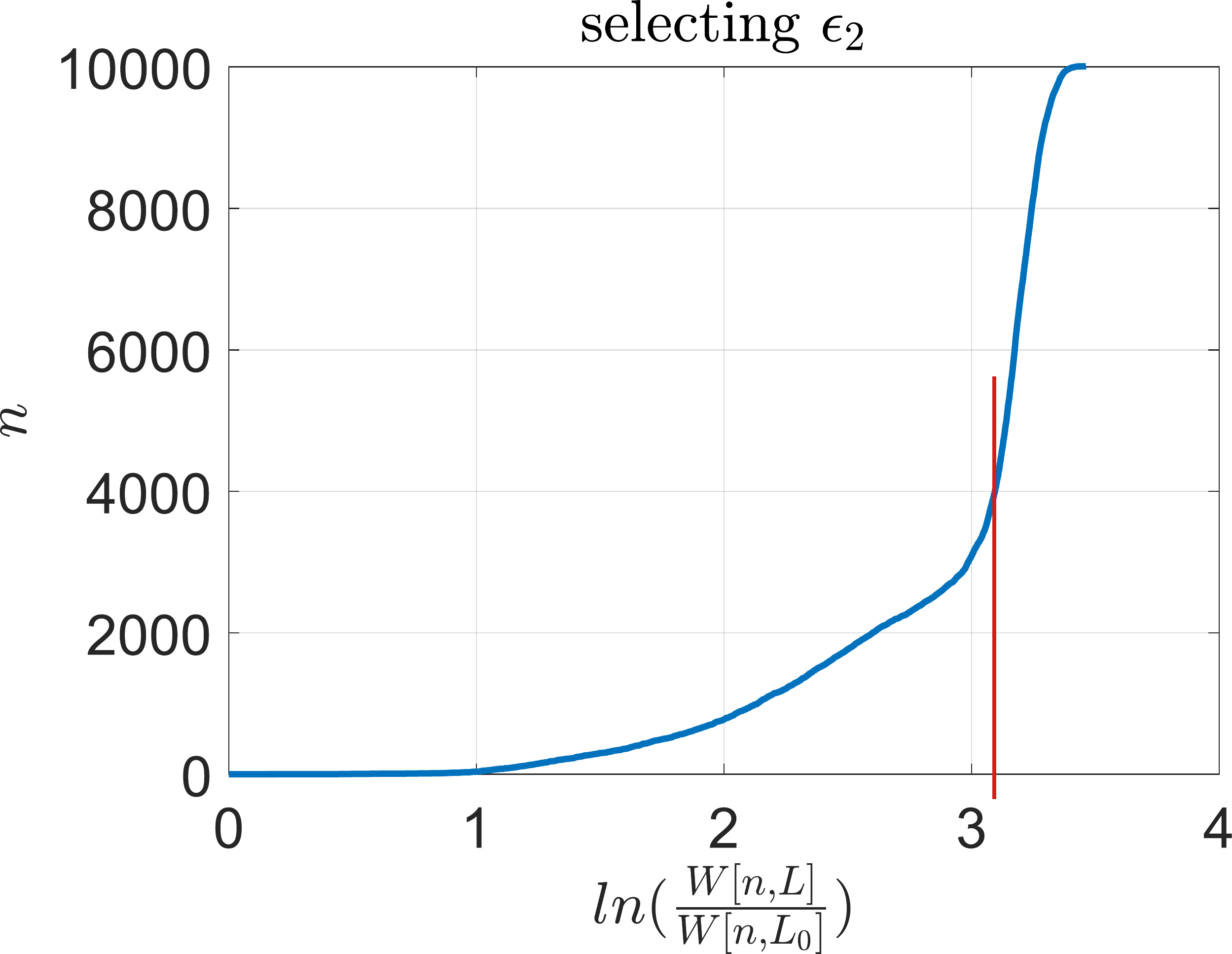}
		\\
		\includegraphics[width=.31\linewidth]{\figs 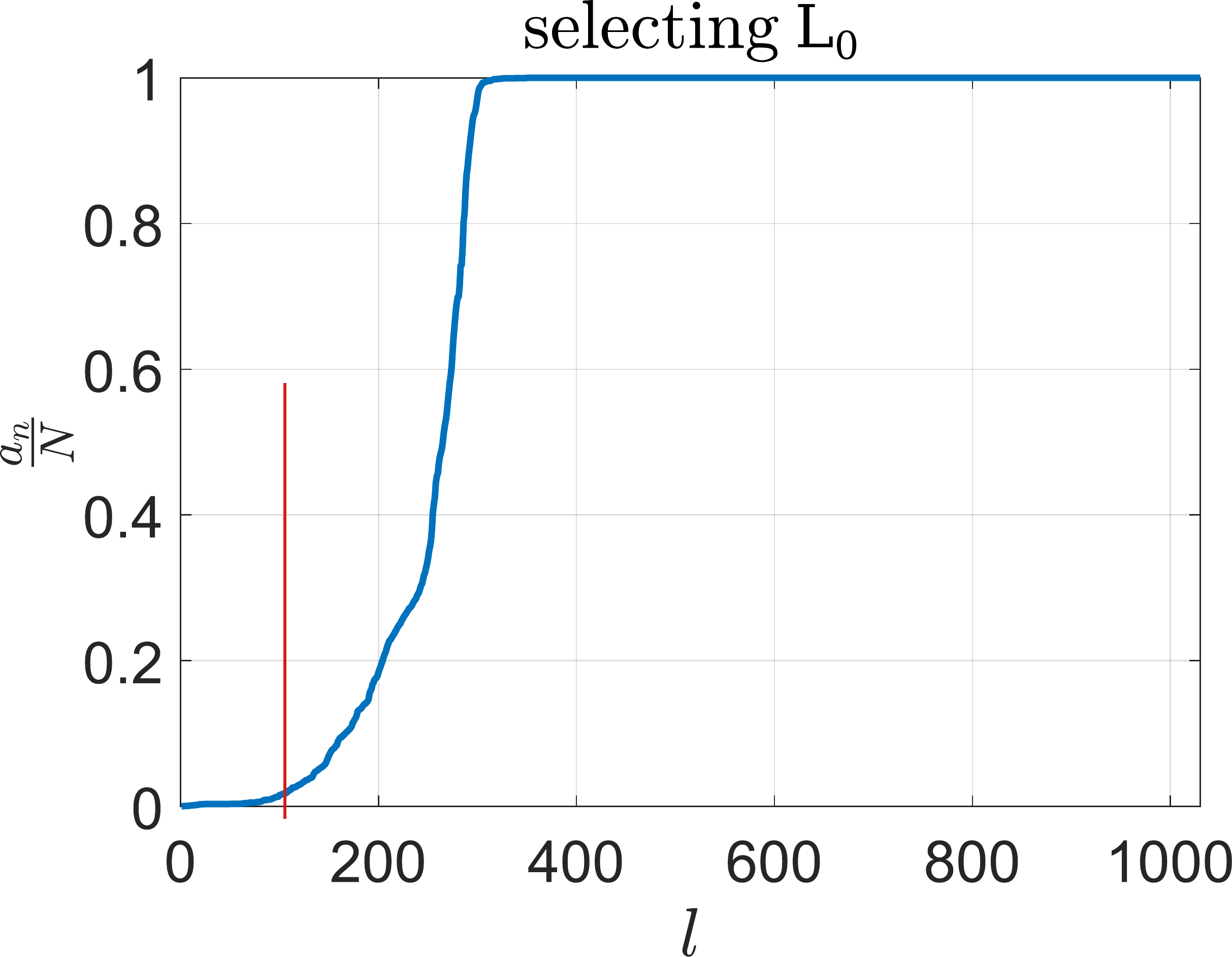}
		\includegraphics[width=.31\linewidth]{\figs 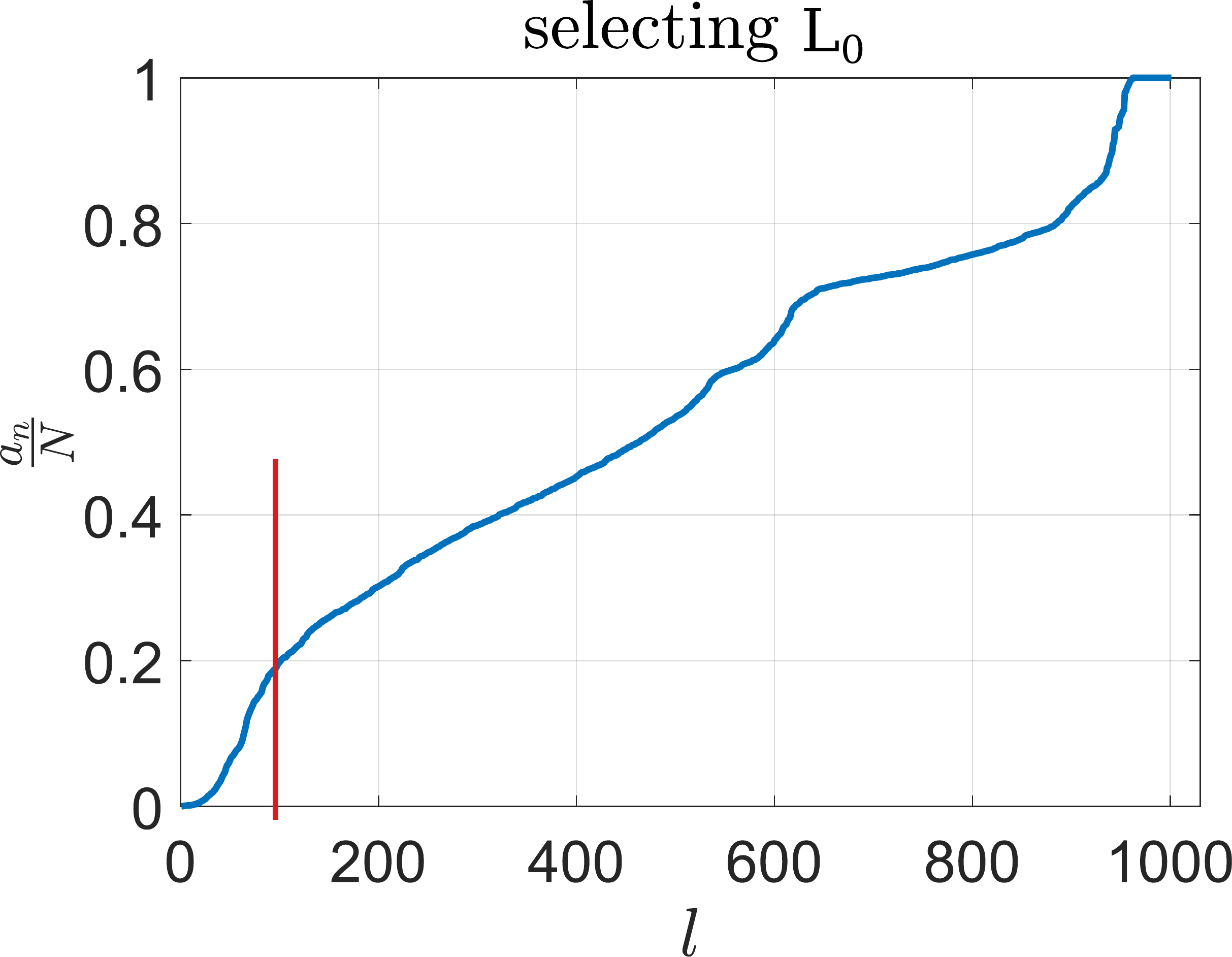}
		\\
		\includegraphics[width=.31\linewidth]{\figs 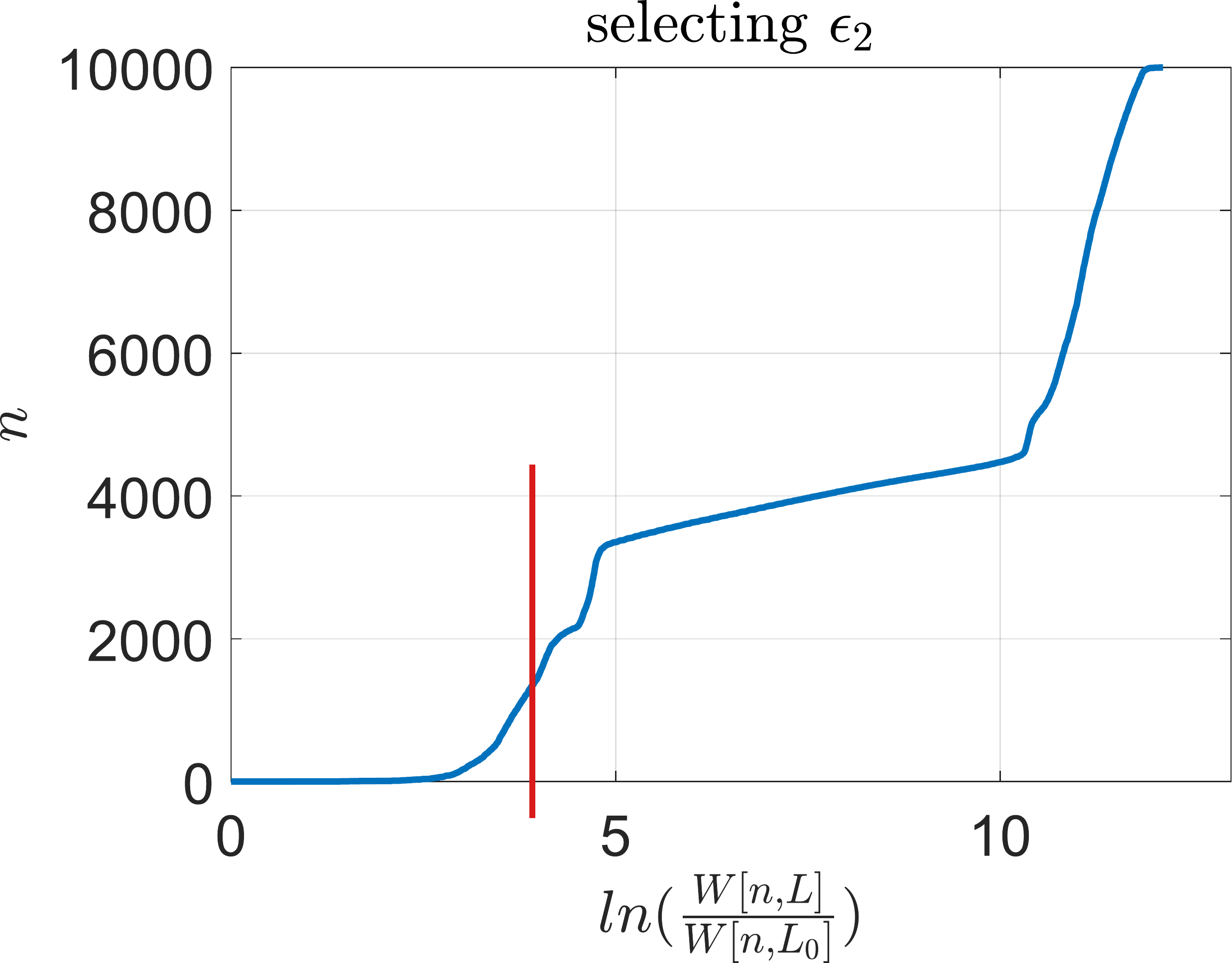}
		\includegraphics[width=.31\linewidth]{\figs 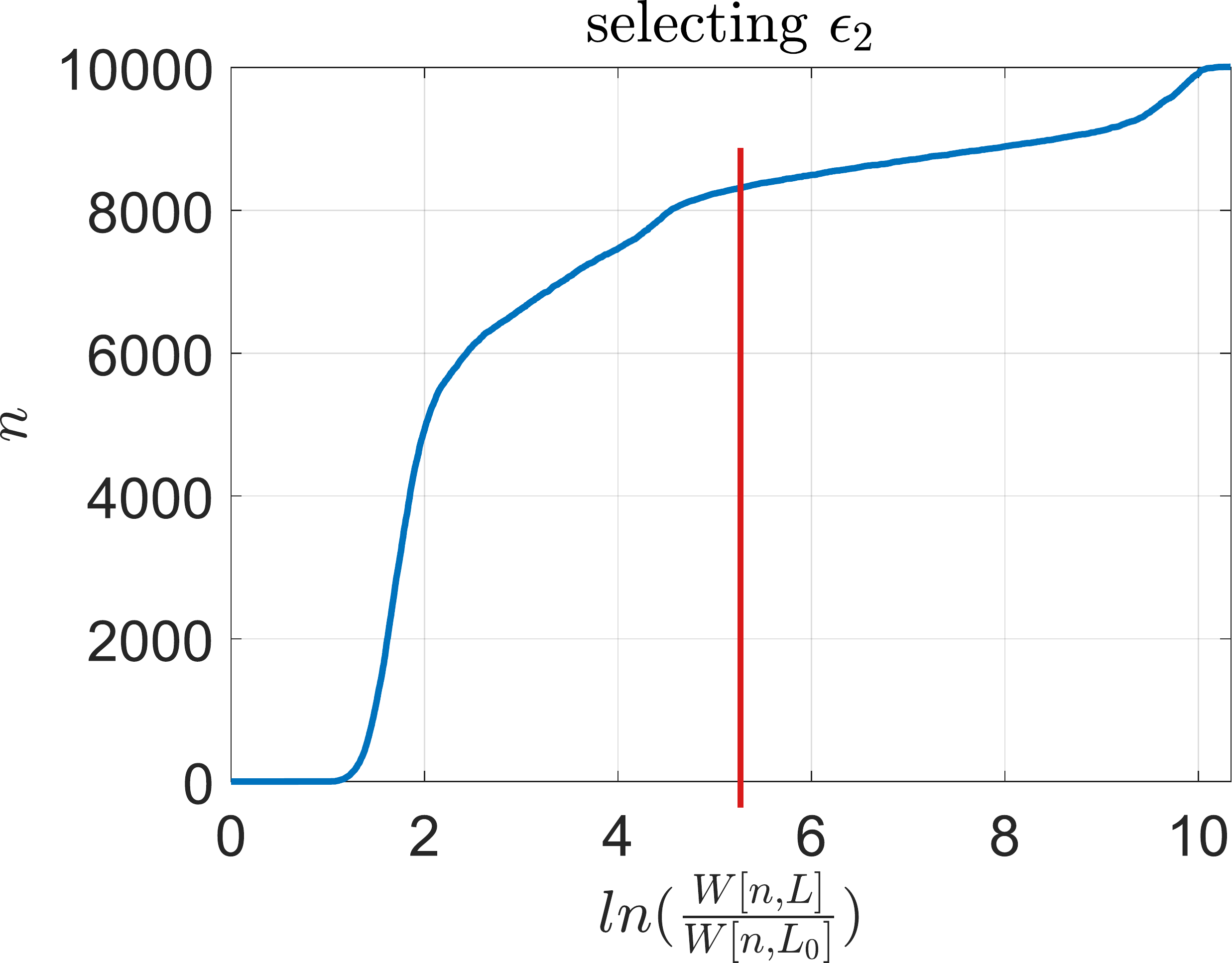}
		\caption{ Choosing the thresholds $L_0, \epsilon_2$ for the three experiments. The left panels correspond to signalized intersections (Alafaya corridor), the middle panels to the freeway (CFX SR 408), and the rightmost panels to the atrial fibrillation data. The red vertical lines correspond to major changes in the slop of the graph and indicate the choice of these parameters. The top panels correspond to the choice of $L_0$, and the bottom panels to $\epsilon_2$.  The lower the value of $\epsilon_2$, the more frequencies get filtered out and the higher the probability of the identification being correct.  See Table~\ref{tab:param} for a complete list of values for $L_0, \epsilon_2$ as well as other parameters. }
		\label{fig:param_selection}
	\end{figure}
	
	We next discuss some important aspects of our numerical reconstruction, brought to light by these results.  
	
	\paragraph{Aliasing} Aliasing is a phenomenon created due to a discrete-time sampling of any continuous time dynamical system. The effect of aliasing is not increased or decreased by the choice of method for analyzing the spectrum. Its effect is seen unavoidably in any numerical analysis of the spectrum. Henceforth, we shall call Koopman eigenfrequencies or eigenvalues as just eigenfrequencies or eigenvalues. 
	
	Recall that an eigenfrequency of $\omega$ of the continuous time dynamical system (or flow) $\SetDef{ \Phi^t : M\to M }{ t\in\real }$ results in the eigenvalue $e^{\iota \Delta \omega}$ for the discrete time analog $\Phi^{\Delta t} : M\to M$. The problem arises due to the $2\pi$-periodicity of the function $x\mapsto e^{\iota x}$. The infinite real line $\real$ is stretched by $1/\Delta t$ and  wrapped around the circle represented by $[0, 2\pi]$. As a result, if there were two eigenfrequencies of the continuous time system $\omega', \omega''$ such that $\omega'' = \omega' + n 2\pi/\Delta t$, both lead to the same eigenvalue :
	\[ e^{\iota \Delta t \omega'} = e^{\iota \Delta t \omega''} = \lambda . \]
	If $\zeta'$ and $\zeta''$ were eigenvectors corresponding to $\omega', \omega''$, then note that
	\[ \zeta' \paran{ \Phi^{\Delta t} x } = \lambda \zeta'(x), \quad \zeta'' \paran{ \Phi^{\Delta t} x } = \lambda \zeta''(x) \]
	Conversely, if $\lambda = e^{\iota \omega}$ is an eigenvalue of the map $\Phi^{\Delta t}$, then it could correspond to any of the infinite set of frequencies $\SetDef{ \omega + n \frac{2\pi}{\Delta t} }{ n = \ldots, -2, -1, 0, 1, 2, \ldots }$. This ambiguity is the source of the problem of aliasing. The angular part $\omega$ of the eigenvalue $\lambda$ is the only representative of an infinite number of translations of a single frequency. If one fixes a fundamental interval such as $[0, 2\pi / \Delta t$ or equivalently $\left[ -\frac{\pi}{\Delta t}, \frac{\pi}{\Delta t} \right]$, and if the flow $\Phi^t$ does not have any eigenfrequencies outside the fundamental interval, then there would not be any effect of aliasing. Thus however is never the case. Remedies to the effect of aliasing in a dynamical context would be
	\begin{enumerate}
		\item Using two different sampling frequencies, in a technique outlined in \citep[][Corr 2]{DasGiannakis_RKHS_2018}.
		\item Passing the signal through a low pass filter that still preserves the dynamics information contained in the signal.
	\end{enumerate}
	We are currently investigating ways to incorporate these techniques into our analysis. One effect of aliasing that is also visible in our calculations is known as \emph{Nystrom phenomenon}. It is seen as the presence of the frequency lowest non-zero $\frac{1}{N \Delta t}$ in our RKHS-based analysis of the signal.
	
	\paragraph{Choice of parameters} The numerical procedure has several parameters associated with it, as summarized in Table \ref{tab:param}. For making the experiments more comparable, we chose the same size $N$ for the data, and the parameter $\epsilon_1$. The bandwidth $\epsilon$ of the kernel and the number of delays $Q$ need to be chosen suited to the data. A small $\epsilon$ captures more geometric information but could lead to a problem of undersampling.  An $\epsilon \approx 0.01k$ seems to be a good choice for most cases. The number of delays $Q$ is set using the correlation-based method \cite{Aguirre1995_delays}. Figure \ref{fig:param_selection} describes a heuristic procedure for selecting the parameters $L_0, \epsilon_1, \epsilon_2$. These three parameters filter the candidate frequencies based on two criteria - their strengths in the original signal as well as the RKHS regularity of their associated waveforms \citep[see][Sec 8]{DasGiannakis_RKHS_2018}.
	
	\paragraph{Learning non-smooth dynamics} Out of the three case studies, the reconstruction of the heart atrial data has the highest error, as seen in Figure~\ref{fig:atrial}. The electrical signals in the heart are intermittent, their spiking or firing behavior makes it a highly non-smooth system. Non-smooth systems pose a major challenge for learning problems. The authors are currently investigating an extension of the model in \eqref{eqn:def:quasi_driven_II} that incorporates this spiking behavior structurally.
	
	\paragraph{Conclusions} We have thus shown the following :
	\begin{enumerate} [(i)]
		\item The use of two different thresholds $\epsilon_1$ and $\epsilon_2$ in our core  Algorithm~\ref{algo:RKHS_kernel} is based on the asymptotic behavior in two different directions, provide a surer guarantee of identification of true eigenfrequencies and the discarding of \emph{spurious} or \emph{pseudo}-spectrum \citep[e.g.][Sec 4.2]{DGJ_compactV_2018}.
		\item Applicability to chaotic dynamics : A unique aspect of our methods is that the methods are applicable to the analysis of signals in which the periodic component is either non-dominant or even absent. For such signals, both DMD-based techniques and Fourier techniques fail to identify the true eigenfrequencies.
		\item Smoothness of reconstruction: an inherent advantage of kernel-based techniques is the easy extrapolation from data to the entire data-space. Moreover, these interpolated functions have the same degree of smoothness as the kernel. The oscillatory behavior of the kernel eigenfunctions increases with the index $l$ and the choice of the spectral resolution parameter controls the smoothness of our interpolation.
		\item Out of sample evaluation : is essentially the task of extrapolation, i.e., evaluating functions reconstructed using the kernel, at points not in the original data-set. At each step of the iteration of \eqref{eqn:reconstruct}, we perform these evaluations using Algorithm~\ref{algo:oss}.
		\item Boundedness of reconstructed dynamics : A crucial advantage that kernel based methods offer over methods such as linear or polynomial regression, is that the interpolation is bounded, due to the decaying nature of the kernel \eqref{eqn:def:Gauss_ker}. All reconstructed dynamical models have some difference with the true system. This difference / defect is inevitable in a learning problem. If the reconstruction of $g_{per}$ is bounded, then it would guarantee that the dynamics under \eqref{eqn:def:quasi_driven_II} would remain bounded, and the deviation of the trajectories also remain bounded.
	\end{enumerate}
	
	We have compared our kernel based method with other spectral estimation techniques in Table~\ref{tab:compare}. One key aspect of using Gaussian kernels is that it leads to the creation of sparse $N\times N$ matrices, which lead to efficient computation and more economical memory usage. The out of sample reconstruction via \eqref{eqn:reconstruct} also gives an explicit formula for the learnt / interpolated function. Another advantage our method derives from the theoretical results of \cite{DasGiannakis_RKHS_2018} is that it can handle generated not only by periodic sources but by systems with purely chaotic or mixed spectrum.
	
	\begin{table}
		\caption{Various techniques for identifying the natural / Koopman eigenfrequencies of dynamical systems, and their performance based on various parameters. }
		\begin{tabularx}{\linewidth}{|L|L|L|L|L|L|L|}
			\hline
			\ & Fourier averaging & EDMD & HDMD & Neural networks with memory & RKHS &\ Non-parametric regression \\ \hline
			Related works & \cite{LangeEtAl2021} \cite{DasJim2017_SuperC} & \cite{WilliamsEtAl15} \cite{SchmidSesterhenn08, Kawahara2016, KutzEtAl16} & \cite{KordaMezic2018}  & \cite{LecunEtAl2015} \cite{YeungEtAl2019}, \cite{HarlimEtAl2021, MaEtAl_2018, Maulik_EtAl_2020, RahmanHasan2020} & \cite{BerryHarlim2017, BerryEtAl2015} \cite{AlxndrGian2020, DasGiannakis_RKHS_2018, DasDimitris_CascadeRKHS_2019} & \cite{Lin2004statistical, HallReimannRice2000, Silverman1984, TompkinsRamos2020} \\ \hline
			Avoids dense $N\times N$ matrix & Y & N & N & Y & Y & Y \\ \hline
			Applicable to chaotic systems & N & N & N & N & Y & N \\ \hline
			Applicable to systems with mixed spectrum & N & N & N & Y & Y & N \\ \hline
			Explicit reconstruction & Y & N & N & N & Y & Y \\ \hline
			Higher accuracy for quasiperiodic systems & Y & Y & Y & N & Y & N \\ \hline
			Low cost of iteration & Y & Y & Y & Y & N & N \\ \hline
		\end{tabularx}
		\label{tab:compare}
	\end{table}
	
	\bibliographystyle{unsrt_inline_url}
	\bibliography{References,SIADS_bib}
\end{document}